\documentclass[fontsize=12pt, DIV=15, parskip=half-]{scrartcl}
\usepackage[utf8]{inputenc}
\usepackage[english]{babel}
\usepackage{tikz}
\usepackage[T1]{fontenc}
\usepackage[lf]{libertine}
\usepackage{textcomp}
\usepackage{microtype}
\usepackage{array}
\usepackage[inline]{enumitem}
\usepackage{amsmath,bm}
\usepackage{amssymb}
\usepackage{amsfonts}
\usepackage[amsmath, thmmarks]{ntheorem}
\usepackage[italic, defaultmathsizes, basic]{mathastext}
\usepackage{stmaryrd}
\usepackage{booktabs}
\usepackage{nicefrac}
\usepackage{algorithm}
\usepackage{algorithmic}
\usepackage[margin=2cm]{geometry}
\usepackage{upgreek}
\usepackage{marvosym}
\usepackage{hyperref}
\hypersetup{
    colorlinks=true,
    linkcolor=blue,
    filecolor=magenta,
    urlcolor=cyan,
}
\usepackage{caption}
\usepackage{subcaption}
\usepackage{array}
\usepackage{lipsum}
\usepackage{xspace}

\usepackage[backend=bibtex,style=numeric,giveninits=true]{biblatex}
\newtheorem{thm}{Theorem}[section]

\newtheorem{lem}[thm]{Lemma}
\newtheorem{prop}[thm]{Proposition}

\newtheorem{rem}[thm]{Remark}
\newtheorem{exa}[thm]{Example}

\newenvironment{prof}[1][{Proof}]{\textbf{#1:} }{\ \rule{0.4em}{0.4em}}

\numberwithin{equation}{section}
  {\pushQED{\qed}\exa}
  {\popQED\endexa}

\newcommand{\R}{\mathbb{R}}
\newcommand{\N}{\mathbb{N}}

\newcommand{\conv}{\textrm{conv }}

\newcommand{\sgn}{\textrm{sgn}}

\newcommand{\Argmin}{\textrm{Argmin }}
\newcommand{\tmenge}[2]{\smash{\{#1\mathbin|#2\}}}
\newcommand{\tnorm}[1]{\mathopen\Vert\smash{#1}\mathclose\Vert}
\newcommand{\eee}{\mathrm{e}}
\newcommand{\0}{\mathbf{0}}
\newcommand{\A}{\mathbf{A}}
\newcommand{\B}{\mathbf{B}}
\newcommand{\vv}{\mathbf{v}}
\newcommand{\I}{\mathbf{I}}
\newcommand{\D}{\mathbf{D}}
\newcommand{\M}{\mathbf{M}}
\newcommand{\NN}{\mathbf{N}}
\newcommand{\C}{\mathbf{C}}
\newcommand{\dd}{\mathbf{d}}
\newcommand{\cc}{\mathbf{c}}
\newcommand{\ttt}{\mathbf{t}}
\newcommand{\yy}{\mathbf{y}}
\newcommand{\xx}{\mathbf{x}}
\newcommand{\zz}{\mathbf{z}}
\newcommand{\hh}{\mathbf{h}}
\newcommand{\rr}{\mathbf{r}}
\newcommand{\kk}{\mathbf{k}}
\newcommand{\sss}{\mathbf{s}}
\newcommand{\bb}{\mathbf{b}}
\newcommand{\uu}{\mathbf{u}}
\newcommand{\pp}{\mathbf{p}}
\newcommand{\aeq}{\Leftrightarrow}

\usepackage{longtable}
\usepackage{lipsum} 
\allowdisplaybreaks
\begin{document}
\title{Newton-type method for bilevel programs with linear lower level problem and application to toll optimization}
\date{}
\author{{\textsf{Floriane Mefo Kue}}\footnote{Department Mathematik, Universit\"{a}t Siegen, Walter-Flex-Str. 3, D-57068 Siegen, Germany \newline   (email: \href{floriane.mkue@uni-siegen.de}{floriane.mkue@uni-siegen.de})},\;
{\textsf{Thorsten Raasch}}\footnote{Department Mathematik, Universit\"{a}t Siegen, Walter-Flex-Str. 3, D-57068 Siegen, Germany \newline (email: \href{raasch@mathematik.uni-siegen.de}{raasch@mathematik.uni-siegen.de})},\;
{\textsf{and Alain B. Zemkoho}}\footnote{School of Mathematical Sciences, University of Southampton, SO17 1BJ Southampton, United Kingdom \newline (email: \href{a.b.zemkoho@soton.ac.uk}{a.b.zemkoho@soton.ac.uk})}}

\maketitle

\begin{abstract}
\noindent \textbf{Abstract.} We consider a bilevel program  involving a linear lower level problem with left-hand-side perturbation. We then consider the Karush-Kuhn-Tucker reformulation of the problem and subsequently build a tractable optimization problem with linear constraints by means of a partial exact penalization. A semismooth system of equations is then generated from the later problem and a Newton-type method is developed to solve it. Finally, we illustrate the convergence and practical implementation of the algorithm on the optimal toll-setting problem in transportation networks. \\

\end{abstract}

\section{Introduction}\label{intro}
Bilevel programming problems are mathematical problems with a special constraint, which is implicitly determined by another optimization problem.
This latter problem, called the follower's or lower level problem, is defined by
\[
\underset{\yy}\min\tmenge{f(\xx,\yy)}{g(\xx,\yy)\leq0},
\]
 where $f: \mathbb{R}^{n}\times\mathbb{R}^{m}\rightarrow\mathbb{R}$ and $g: \mathbb{R}^{n}\times\mathbb{R}^{m}\rightarrow\mathbb{R}^{l}$ represent the lower level objective and constraint functions, respectively.
 The bilevel programming problem can be formally described as 
 \begin{equation} \label{Prelim_Pb}
\begin{array}{rl}
   \underset{\xx,\yy}\min & F(\xx,\yy)\\
 \mbox{ s.t. }& G(\xx)\leq0,\quad \yy \in \Psi(\xx):= \underset{\yy}{\Argmin} \tmenge{f(\xx,\yy)}{g(\xx,\yy)\leq0},
\end{array}
\end{equation}
where, similarly, the  functions $F: \mathbb{R}^{n}\times\mathbb{R}^{m}\rightarrow\mathbb{R}$ and $G: \mathbb{R}^{n}\rightarrow\mathbb{R}^{k}$ correspond to the upper level objective and constraint functions, respectively.
Problem \eqref{Prelim_Pb} is called the upper level or leader's problem. It is often referred
to as the optimistic formulation of the bilevel programming problem \cite{dempe2002foundations}. The variables $\xx$ and $\yy$ are the upper and lower level variables or upper and lower level decisions, respectively.
In this order of ideas, our focus in this paper will be on the simple bilevel program
 \begin{equation} \label{investigated_Problem}
\begin{array}{rl}
   \underset{\xx,\yy}\min & F(\xx,\yy)\\
 \mbox{ s.t. }& \D\xx\leq \dd,\quad \yy \in \Psi(\xx),
\end{array}
\end{equation}
with the lower level optimal solution set-valued mapping $\Psi$ defined by
\begin{equation} \label{lower-level_problem}
\Psi(\xx):= \underset{\yy}{\Argmin}\tmenge{\xx^\top \yy}{\A \yy\leq \bb},
\end{equation}
where $\xx,\yy\in \R^n$, $\dd\in \R^m$, $\bb\in \R^l$ and all the remaining vectors and  matrices  are of appropriate dimensions. Note that the matrix $\A$ involved in \eqref{lower-level_problem} does not depend on the parameter $\xx$ and the function $F$ is assumed to be twice continuously differentiable. The term  \emph{simple bilevel program} for this class of problem was coined in the paper \cite{lin2014solving}, referring to the fact that the upper (resp. lower) level feasible set is independent from $\yy$ (resp. $\xx$). However, it is important to mention that the same expression was initially used in \cite{dempeSIMPLE2010optimality} to name a completely different class of problem, that can be traced back to \cite{mangasarian1979nonlinear} and has been investigated in various publications; see, e.g., \cite{cabot2005proximal, ferris1991finite, solodov2007explicit, shehu2019inertial}.

It is clear that the lower level problem \eqref{lower-level_problem} is linear (w.r.t. $\yy$), but with a left-hand-side perturbation. However, problem \eqref{investigated_Problem} is still a significantly difficult class of problem; to have a taste of this, consider the following example, where $F$ is any real-valued function and $\xx$ and $\yy$ are one-dimensional:
\begin{equation} \label{exa}
\begin{array}{rl}
   \underset{\xx, \yy}\min & F(\xx,\yy)\\
 \mbox{ s.t. }& \xx\in[-1,1],\quad \yy \in \Argmin\{\xx \yy\,|\;\yy\in[0,1]\}.
\end{array}
\end{equation}
The feasible set of the problem is depicted in the following picture:
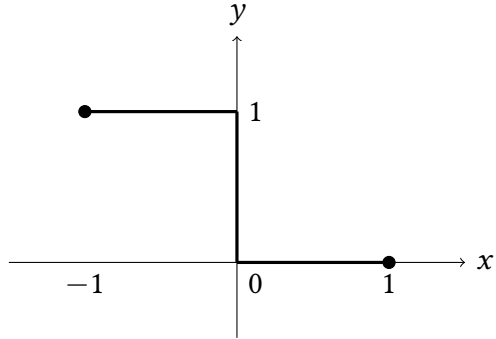
\begin{center}
\begin{tikzpicture}
\draw[->] (-3,0) -- (3,0);
\draw (3,0) node[right] {$x$};
\draw (2,0) node[below] {$1$};
\draw [->] (0,-1) -- (0,3);
\draw[very thick] (0,0) -- (2,0);
\draw (0,3) node[above] {$y$};
\draw (0,2) node[right] {$1$};
\draw (-2,0) node[below] {$-1$};
\draw (0,0) node[below right] {$0$};
\draw[very thick] (0,0) -- (0,2);
\draw[very thick] (0,2) -- (-2,2);
\draw [fill] (2,0) circle (0.08cm);
\draw [fill] (-2,2) circle (0.08cm);
\end{tikzpicture}
\end{center}
Clearly, the feasible set of problem \eqref{exa}, represented by the thick line segments, is the union of convex polyhedral sets; but it is not a convex set itself. Hence, finding a global optimal solution for problem \eqref{investigated_Problem}--\eqref{lower-level_problem} is generally not an easy task. This is why our focus in this paper will be on computing stationary points of the problem, which can potentially be locally optimal.

Note that if we replace the lower level problem \eqref{lower-level_problem} with the following one, where the perturbation is instead in the right-hand-side, we get a much easier optimization problem:
 \begin{equation} \label{lower-level_problem1}
\Psi(\xx):= \underset{\yy}{\Argmin}\tmenge{\cc^\top \yy}{\A \yy\leq \xx}.
\end{equation}
Problem \eqref{lower-level_problem1} is \emph{fully convex}, in the sense the objective and constraint functions are convex in $(\xx, \yy)$. Thanks to this,  some specific algorithms work well for \eqref{investigated_Problem}, \eqref{lower-level_problem1}, but cannot be implemented on \eqref{investigated_Problem}, \eqref{lower-level_problem}; cf. \cite[Chapter 3]{dempe2002foundations}. Also, the optimal value function of problem \eqref{lower-level_problem1} is convex and \eqref{investigated_Problem}, \eqref{lower-level_problem1} is partially calm  \cite{Mehlitz2020ANO}. The latter properties are critical in addressing problem \eqref{investigated_Problem}, \eqref{lower-level_problem1} theoretically and numerically. But unfortunately, they cannot hold for \eqref{investigated_Problem}, \eqref{lower-level_problem}. This probably justifies why problem \eqref{investigated_Problem}, \eqref{lower-level_problem1} has been at the center of attention in the development of solution algorithms for linear bilevel programs; see, e.g., \cite{CalveteGale2020BilevelBook} for the most recent literature review on the subject.
Note that the approach to be developed in this paper can be applied to problem \eqref{investigated_Problem}, \eqref{lower-level_problem1} as well.

There are various practical bilevel programs with linear lower-level problems that are not covered by the model \eqref{investigated_Problem}, \eqref{lower-level_problem1}. Among them, the most prominent is probably the toll optimization problem, also known as optimal toll-setting or network pricing problem \cite{brotcorne2000bilevel,dewez2008new,labbe1998bilevel} that we are going to use in Section \ref{sect. app} to illustrate the practical implementation of the algorithm developed in this paper.

To tackle problem \eqref{investigated_Problem}, \eqref{lower-level_problem}, we consider its Karush-Kuhn-Tucker (KKT) reformulation, as the required lower level convexity and regularity are both automatically satisfied \cite{dempe2012bilevel}. Then, using a certain tractable transformation process, we construct a partial exact penalization of problem \eqref{investigated_Problem}, \eqref{lower-level_problem} and show  its close link with the KKT reformulation; cf. Section \ref{Single-level reformulation}. In the context of this penalized problem, we investigate the semismooth Newton method (see, e.g., \cite{facchinei1997new, nesterov2013introductory, qi1993nonsmooth}) in Section \ref{Semismooth Newton method}, where sufficient conditions ensuring its convergence are established and illustrated on some examples. The method is then implemented on the toll optimization problem in Section \ref{sect. app}.

It is important to recall that semismooth Newton-type methods based on the KKT or value function reformulations of general versions of the bilevel optimization problem \eqref{investigated_Problem}, \eqref{lower-level_problem} have been studied  recently in \cite{fischer2019semismooth, fliege2020gauss, zemkoho2020theoretical}. These works typically require the partial calmness property, which does not necessarily hold for problem \eqref{investigated_Problem}, \eqref{lower-level_problem}. 
Moreover, for this problem specifically, very little work tailored to it has been done in the literature. We are aware of the paper \cite{dempe2014solution}, where an algorithm based on the optimal value reformulation to compute local and global optimal solutions is derived. In \cite{dempe2017solving} the authors propose an approach to solve \eqref{investigated_Problem}, \eqref{lower-level_problem} in the discrete case through the optimal value reformulation approach as well. For a more general version of the problem, the authors in \cite{lin2014solving}  design a smoothing projected gradient algorithm for a simple bilevel programs with a nonconvex lower level program using the optimal value reformulation approach.

Before, we go deep in the analysis as described above, we start with some preliminary elements in the next section.

\section{Preliminaries}
Here, we introduce some preliminary concepts  about the semismooth Newton method that will be used later in this paper. First, note that a function $\Phi:\R^n\rightarrow \mathbb{R}^{m}$ is said to be locally Lipschitz continuous around $\xx^{0}\in \R^n$ if there exist $\alpha>0$ and $L>0$ such that the following condition holds:
\[
\forall \xx,\xx'\in \xx^0+\alpha\mathbb{B}_{\R^n}\colon\quad\Vert \Phi(\xx)-\Phi(\xx')\Vert_{\R^m}\leq L\Vert \xx-\xx'\Vert_{\R^n}.
\]
The number $L$ is called the Lipschitz constant. $\Phi$ will be said to be locally
Lipschitz continuous if it is locally Lipschitz continuous around every point of $\R^n$. It is said to be Lipschitz continuous if the above inequality holds with $\alpha=\infty$. Recall that any convex function is locally Lipschitz continuous on the relative interior of its domain.

Consider a locally Lipschitz continuous function $\Phi:\R^n\rightarrow \mathbb{R}^{m}$. According to Radamacher's theorem \cite{clarke1990optimization}, $\Phi$ is differentiable almost everywhere. Let $D_\Phi$ denote the set of all the points where $\Phi$ is differentiable.  Then, we can define the Bouligand subdifferential of $\Phi$ at $\xx\in\R^n$ by
\[
\partial_B\Phi(\xx):=\tmenge{\C\in\R^{n\times m}}{\C=\lim_{k\to\infty}\nabla\Phi(\xx_k), \;\xx_k\to \xx,\;\xx_k\in D_\Phi}.
\]
The generalised Jacobian of $\Phi$ in the sense of Clarke \cite{clarke1990optimization} is given by
\[
\partial\Phi(\xx):=\conv\partial_B\Phi(\xx),
\]
where \emph{conv} stands for the convex hull.
The postulated local Lipschitz continuity property of $\Phi$ around $\xx$ guarantees that the set $\partial\Phi( \xx)$ is nonempty and compact.
If $\Phi$ is a convex function, then $\partial\Phi(\xx)$ coincides with the subdifferential in the sense of convex analysis.

Next we introduce the notion of semismoothness introduced in \cite{mifflin1977semismooth} and extented  to vector-valued functions in \cite{qi1993convergence, qi1993nonsmooth}.	
 A function $\Phi:\R^n\rightarrow \mathbb{R}^{m}$ is said to be semismooth at a given point $\xx$ if $\Phi$ is directional differentiable at $\xx$ and for any $\C\in\partial \Phi(\xx+\hh),\hh\to0$
 \[
 \Phi(\xx+\hh)-\Phi(\xx)-\C\hh=o(\tnorm{\hh}).
 \]
 The function $\Phi$ is said to be strongly  semismooth  at $\xx$ if $\Phi$ is semismooth at $\xx$ and it holds that
 \[
 \Phi(\xx+\hh)-\Phi(\xx)-\C\hh=O(\tnorm{\hh}^2),\;\forall\C\in\partial \Phi(\xx+\hh),\hh\to0.
 \]
 Any continuously differentiable is obviously semismooth. But this is not necessarily the case for a locally Lipschitz continuous  \cite{mifflin1977semismooth}. However, piecewise linear functions are strongly semismooth and a special case from such a class, that will be of a particular interest in this paper is 
 $t\mapsto\max(0, t)$. Its Clarke subdifferential can be obtained as
 \[
 v\in \partial \max(0,t)\aeq v\in\begin{cases}
 \{1\}&\;\text{if} \;\;\; t>0,\\
 \{0\}&\;\text{if} \;\;\; t<0,\\
 [0,1]&\;\text{if} \;\;\; t=0.\\
 \end{cases}
 \]
\indent Throughout the paper, we use $\mathbb{R}^{n}_{+}$ and  $\mathbb{R}^{n\times m}$ to denote the cone of  $n$-dimensional real-valued vectors with nonnegative components and the  set  of  real-valued  matrices  with $n$ rows  and $m$ columns, respectively. Furthermore we denote by $\0$, $\I$, and $\eee$ the null matrix, the identity matrix and the vector of ones of appropriate dimension, respectively.

\section{Single-level reformulation}\label{Single-level reformulation}
Considering problem \eqref{investigated_Problem}, \eqref{lower-level_problem}, it is obvious that its lower level problem is both convex and linear w.r.t. the lower level variable $\yy$. Hence, the KKT reformulation of the problem can be written as follows without any additional requirement:
 \begin{equation}\label{KKT_Reformulation}
\begin{array}{rl}
   \underset{\xx, \yy, \zz}\min & F(\xx,\yy)\\
 \mbox{ s.t. }& \D\xx\leq \dd, \;\; \A^\top \zz+\xx=0,\\
&\A\yy\leq \bb, \;\; \zz\geq 0, \;\; \zz^\top(\A\yy-\bb)=0.
\end{array}
\end{equation}
Reference \cite{dempe2012bilevel} provides a detailed analysis of the relationship between this reformulation and the original problem \eqref{investigated_Problem}, \eqref{lower-level_problem}. One of the main issues in solving problem \eqref{KKT_Reformulation} with standard continuous optimization techniques is the presence of the complementarity conditions
\begin{equation}\label{Complementarity}
  \left(\bb-\A\yy\right)_i\geq0, \;\;\zz_i\geq0, \;\; \zz_i\left(\bb-\A\yy\right)_i=0 \;\;\mbox{ for } \;\;i=1, \ldots, l
\end{equation}
in the feasible set, which cause the failure of well-known constraint qualifications (see, e.g., \cite{dempe2012karush}). To deal with this issue here, we start by considering the reformulation
\begin{equation}\label{Complementarity Reformulation}
  \underset{(\rr_i,\sss_i)\in T_i}\min \rr_i\zz_i+\sss_i(\bb-\A\yy)_i=0 \;\;\mbox{ for } \;\;i=1, \ldots, l,
\end{equation}
of the complementarity conditions \eqref{Complementarity}, where, for $i=1, \ldots, l$, the set $T_i$ is defined by
\begin{equation}\label{Set T}
 T_i:=\tmenge{(a,b)\in\R^{2}}{a\geq0, \;\, b\geq0, \;\, a+b=1}.
\end{equation}
Note that transformation \eqref{Complementarity Reformulation} obviously follows from the fact \eqref{Complementarity} is equivalent to considering the fact that \eqref{Complementarity} is identical to $\min~(\zz_i,\, (\bb-\A\yy)_i)=0$ for $i=1, \ldots, l$.


Based on \eqref{Complementarity Reformulation}, we consider the following penalization of \eqref{KKT_Reformulation}:
 \begin{equation}\label{Penalized_Pb}\tag{P$_\alpha$}
\begin{array}{rl}
   \underset{\xx,\yy,\zz}\min & F(\xx,\yy)+\alpha\pi(\yy,\zz)\\
 \mbox{ s.t. }& (\xx,\yy,\zz)\in Z_1:=\tmenge{(\xx,\yy,\zz)}{\D\xx\leq \dd,\;\A\yy\leq \bb,\;\A^\top\zz +\xx=0, \zz\geq 0},
\end{array}
\end{equation}
where $\alpha>0$ represents the penalization parameter and
\[
\pi(\yy,\zz):=\sum_{i=1}^{l}\underset{(\rr_i,\sss_i)\in T_i}\min~\rr_i\zz_i+\sss_i(\bb-\A\yy)_i\geq 0.
\]
To establish the relationship between  \eqref{Penalized_Pb} and \eqref{KKT_Reformulation}, let  $T:=\prod^l_{i=1}T_i$
and let $Z_1^*$  be the set of vertices of the polyhedra $Z_1$. 
We have the following result inspired by \cite[Proposition 3.1]{mangasarian1997exact}. 
\begin{prop}\label{Prop-global-solut-equivalence}
 Assume that  for all $\alpha>0$, problem \eqref{Penalized_Pb} possesses an optimal solution which belongs to $Z_1^*$. Then the following two statements are valid:
\begin{enumerate}
\item[\emph{(i)}] There exists a scalar $\bar\alpha>0$ such that for all $\alpha\geq\bar\alpha$ any global optimal solution  $(\bar \xx_\alpha,\bar \yy_\alpha,\bar \zz_\alpha)$ of problem \eqref{Penalized_Pb} which  belongs to $Z_1^*$ solves problem \eqref{KKT_Reformulation} globally.
\item[\emph{(ii)}]  There exists a scalar $\bar\alpha>0$ such that for all $\alpha\geq\bar\alpha$,  any global  optimal solution $(\bar \xx,\bar \yy,\bar \zz)$ of problem \eqref{KKT_Reformulation}  optimally solves problem \eqref{Penalized_Pb} globally.
\end{enumerate}
\end{prop}
\begin{prof}
For (i), we start by noting that any global solution $(\bar \xx_\alpha,\bar \yy_\alpha,\bar \zz_\alpha)$ of \eqref{Penalized_Pb} which  belongs to $Z_1^*$ and  satisfy $\pi(\bar \yy_\alpha,\bar \zz_\alpha)=0$ solves \eqref{KKT_Reformulation} globally. In fact let $(\xx,\yy,\zz)$ be a feasible point of \eqref{KKT_Reformulation}. We have
\[
F(\xx,\yy)=F(\xx,\yy)+\alpha\pi(\yy,\zz)\geq F(\bar \xx_\alpha,\bar \yy_\alpha)+\alpha\pi(\bar \yy_\alpha,\bar \zz_\alpha)=F(\bar \xx_\alpha,\bar \yy_\alpha).
\]
Hence it suffices to show that for all sufficiently large $\alpha$, every global solution  $(\bar \xx_\alpha,\bar \yy_\alpha,\bar \zz_\alpha)$ of \eqref{Penalized_Pb} which  is an extreme point of $Z_1$ satisfies $\pi(\bar \yy_\alpha,\bar \zz_\alpha)=0.$ We suppose that this is not true. Then there exists sequences $\{\alpha_k\}\to\infty$ and $\{X_k:=(\xx_k,\yy_k,\zz_k)\}_k\subseteq Z_1^*$ where $X_k$ is optimal solution of \eqref{Penalized_Pb} and $\pi(\yy_k,\zz_k)>0$.
Due to the fact that the set of extreme point $Z_1^*$ is finite we have
\[
\beta=\underset{k}\inf~\pi(\yy_k,\zz_k)>0 \;\;\mbox{ and } \;\;\rho=\underset{k}\inf~F(\xx_k,\yy_k)>-\infty.
\]
Then for any feasible point $(\xx,\yy,\zz)$ of \eqref{KKT_Reformulation} it holds that
\[
F(\xx,\yy) \geq F(\xx_k,\yy_k)+\alpha_k\pi(\yy_k,\zz_k) \geq \rho+\alpha_k\beta.
\]
Therefore, for $k\to\infty$ we get $F(\xx,\yy)=\infty$  and this is a contradiction. Consequently there exists a scalar $\bar\alpha>0$ such that for all $\alpha\geq\bar\alpha$ any optimal solution of  \eqref{Penalized_Pb} which belongs to $Z_1^*$ solves \eqref{KKT_Reformulation}.

As for (ii), let $\bar\alpha>0$ be the penalization parameter whose the existence has been proved in the first assertion  of this proof. Let  $(\bar \xx,\bar \yy,\bar \zz)$ be a global optimal solution of \eqref{KKT_Reformulation}. Let $\alpha\geq\bar \alpha$ and $\bar X_{\alpha}:=(\bar\xx_{\alpha},\bar\yy_{\alpha},\bar\zz_{\alpha})$ be a global  optimal solution of \eqref{Penalized_Pb} which belongs to $Z_1^*$ (such point exists from our assumption.). Then $(\bar\xx_{\alpha},\bar\yy_{\alpha},\bar\zz_{\alpha})$ solves \eqref{KKT_Reformulation} from the first-proved assertion.  We want to show  that  $(\bar \xx,\bar \yy,\bar \zz)$ solves \eqref{Penalized_Pb}.
For any feasible solution $X_{\alpha}=(\xx_{\alpha},\yy_{\alpha},\zz_{\alpha})$ of \eqref{Penalized_Pb}, for $\alpha\geq\bar\alpha$,
\begin{align*}
F(\xx_{\alpha},\yy_{\alpha})+\alpha\pi(\yy_{\alpha},\zz_{\alpha})&\geq F(\bar \xx_{\alpha},\bar \yy_{\alpha})+\alpha\pi(\bar\yy_{\alpha},\bar\zz_{\alpha})\\
&=F(\bar\xx_{\alpha},\bar\yy_{\alpha})\\
&=F(\bar \xx,\bar \yy)\\
&=F(\bar \xx,\bar \yy)+\alpha\pi(\bar \yy,\bar \zz),\quad
\end{align*}
where the first and second equality are due to the fact that on the first hand $\bar X_\alpha$ belongs to $Z^*_1$ (see also the proof of the first assertion) and on the second hand  $(\bar\xx_{\alpha},\bar\yy_{\alpha},\bar\zz_{\alpha})$ and $(\bar \xx,\bar \yy,\bar \zz)$ solve \eqref{KKT_Reformulation}.
Therefore the point $(\bar \xx,\bar \yy,\bar \zz)$  solves \eqref{Penalized_Pb}.
\end{prof}

\begin{rem}\label{Remark}
One way to ensure the existence of solution of \eqref{Penalized_Pb} which belongs to $Z^*_1$ is to assume that the function $F$ is  concave, the set $Z_1$ is nonempty and $F$ is bounded below on $Z_1$. In fact
since the set $Z_1$ is a polyhedral containing no lines (it is not difficult to verify it since $Z_1$ is supposed to have at least one extreme point) and that for all $\alpha$, the mapping $(\xx,\yy,\zz)\mapsto F(\xx,\yy)+\alpha\pi(\yy,\zz)$ is concave, we conclude using Corollary 32.3.4 in \cite{rockafellar1970convex} or  \cite[Theorem 3.4.7]{bazaraa2013nonlinear} that \eqref{Penalized_Pb} is well defined  and there exists an optimal solution which belongs to $Z_1^*$.
\end{rem}

\begin{prop}\label{Prop-local-solut-equivalence}
 Assume that  for all $\alpha>0$, problem \eqref{Penalized_Pb} possesses an optimal solution which belongs to $Z_1^*$. Then there exists a scalar $\bar\alpha>0$ such that for all $\alpha\geq\bar\alpha$, the following two assertions hold:
\begin{enumerate}
\item[\emph{(i)}] Any local optimal solution  $(\bar \xx_\alpha,\bar \yy_\alpha,\bar \zz_\alpha)$ of  \eqref{Penalized_Pb} solves \eqref{KKT_Reformulation} locally.
\item[\emph{(ii)}] Suppose that  $(\bar \xx,\bar \yy,\bar \zz)$ is local optimal solution  of problem \eqref{KKT_Reformulation} for all $\bar \zz\in\R^l $ with $\pi(\bar \yy,\bar\zz)=0$. The point $(\bar \xx,\bar \yy,\bar \zz)$  solves \eqref{Penalized_Pb} locally if from any feasible point $(\xx,\yy,\zz)$ of problem \eqref{Penalized_Pb}, we can find a point $\tilde \zz\in \R^{l}$ such that the point $(\xx,\yy,\tilde \zz)$ is feasible for problem \eqref{KKT_Reformulation}.
\end{enumerate}
\end{prop}
\begin{prof}
Let $\bar\alpha>0$ as given in the first part of the proof of Proposition \ref{Prop-global-solut-equivalence}.

For (i), let $(\bar \xx_\alpha,\bar \yy_\alpha,\bar \zz_\alpha)$  be any local solution of \eqref{Penalized_Pb}. Then $(\bar \xx_\alpha,\bar \yy_\alpha,\bar \zz_\alpha)$  is also a local solution of \eqref{KKT_Reformulation}. In fact we suppose that it is not true. Then  there exists a sequence of feasible points $(\bar \xx^n,\bar \yy^n,\bar \zz^n)_n$ for \eqref{KKT_Reformulation} converging to $(\bar \xx_\alpha,\bar \yy_\alpha,\bar \zz_\alpha)$ such that $F(\bar \xx^n,\bar \yy^n)<F(\bar \xx_\alpha,\bar \yy_\alpha)$ for all $n\in \N$.
Since $(\bar \xx^n,\bar \yy^n,\bar \zz^n)$ is feasible for \eqref{KKT_Reformulation}, we get $\pi(\bar \yy^n,\bar \zz^n)=0$. Therefore, we have the following inequalities, which lead to a contradiction:
\[
F(\bar \xx^n,\bar \yy^n)+\alpha\pi(\bar \yy^n,\bar \zz^n)=F(\bar \xx^n,\bar \yy^n)<F(\bar \xx_\alpha,\bar \yy_\alpha)\leq F(\bar \xx_\alpha,\bar \yy_\alpha)+\alpha\pi(\bar \yy_\alpha,\bar \zz_\alpha).
\]

In the case of (ii), let $(\bar \xx,\bar \yy,\bar \zz)$ be a local optimal solution of \eqref{KKT_Reformulation}. We suppose by contradiction that  $(\bar \xx,\bar \yy,\bar \zz)$ is not a local optimal solution of \eqref{Penalized_Pb}. Then there exists a sequence of feasible points $(\bar \xx_\alpha^n,\bar \yy_\alpha^n,\bar \zz_\alpha^n)_n$ of \eqref{Penalized_Pb} converging to $(\bar \xx,\bar \yy,\bar \zz)$ such that
\begin{equation*}
F(\bar \xx_\alpha^n,\bar \yy_\alpha^n)+\alpha\pi(\bar \yy_\alpha^n,\bar \zz_\alpha^n)<F(\bar \xx,\bar \yy)+\alpha\pi(\bar \yy,\bar \zz),\;\forall  n\in \N.
\end{equation*}
  Since $(\bar \xx,\bar \yy,\bar \zz)$ is feasible for \eqref{KKT_Reformulation}, we get $\pi(\bar \yy,\bar \zz)=0$. We can then deduce that \mbox{$F(\bar \xx_\alpha^n,\bar \yy_\alpha^n)<F(\bar \xx,\bar \yy)$} because $\pi(\bar \yy_\alpha^n,\bar \zz_\alpha^n)\geq0$. Moreover, we can  find  for each $n\in\N$ a point $\tilde z_n\in\R^l$  such that the point $(\bar \xx_\alpha^n,\bar \yy_\alpha^n,\tilde \zz_\alpha^n)$ is feasible for \eqref{KKT_Reformulation}.  Due to the fact that the point-to-set mapping
  \[
  (\xx,\yy)\mapsto \left\{\zz\geq0\,|\, \A^{\top}\zz+\xx=\0,\,\zz^{\top}(\A\yy-\bb)=0\right\}
  \]
  is upper semicontinuous \cite{robinson1982generalized}, the sequence $\tilde \zz_n$ has an accumulation point $\tilde \zz$ such that $\pi(\bar \yy,\tilde\zz)=0$. This violates the assumption.
\end{prof}

In order to render \eqref{Penalized_Pb} more tractable we use the following proposition to get rid of the operator min of the function $\pi$ which appears in the objective function of \eqref{Penalized_Pb}.
\begin{prop}\label{Prop-problem-equivalent}
For any $\alpha>0$, problem \eqref{Penalized_Pb} is globally equivalent to
 \begin{equation}\label{Penalized_Pb1} \tag{Q$_{\alpha}$}
\begin{array}{rl}
   \underset{\xx,\yy,\zz,\rr,\sss}\min & F(\xx,\yy)+\alpha(\rr^{\top}\zz+\sss^{\top}(\bb-\A\yy)) \;\; \mbox{ s.t. } \;\; (\xx,\yy,\zz,\rr,\sss)\in Z_1\times T
\end{array}
\end{equation}
and the following assertions hold:
\begin{enumerate}
\item[\emph{(i)}] From any local optimal solution $(\bar \xx_\alpha,\bar \yy_\alpha,\bar \zz_\alpha)$ of problem \eqref{Penalized_Pb}, the point  $(\bar \xx_\alpha,\bar \yy_\alpha,\bar \zz_\alpha,\bar \rr_\alpha,\bar \sss_\alpha)$, with the component $(\bar \rr_\alpha,\bar \sss_\alpha)\in\Argmin\pi(\bar \yy_\alpha,\bar \zz_\alpha)$, is local optimal solution of \eqref{Penalized_Pb1} as well.
\item[\emph{(ii)}]  From any local optimal solution $(\bar \xx_\alpha,\bar \yy_\alpha,\bar \zz_\alpha,\bar \rr_\alpha,\bar \sss_\alpha)$  of problem \eqref{Penalized_Pb1} with $(\bar \rr_\alpha,\bar \sss_\alpha)\in\Argmin\pi( \yy, \zz)$ for all $(\yy, \zz)$ such that $( \bar \xx_\alpha,\yy, \zz)\in Z_1$, the point $(\bar \xx_\alpha,\bar \yy_\alpha,\bar \zz_\alpha)$ is a local solution of problem \eqref{Penalized_Pb}.
\end{enumerate}
\end{prop}
\begin{prof}
The initial statement on the global equivalence between problems \eqref{Penalized_Pb} and \eqref{Penalized_Pb1} is obvious. As for (i), let $(\bar \xx_\alpha,\bar \yy_\alpha,\bar \zz_\alpha)$  be a local  optimal solution of problem \eqref{Penalized_Pb}. We want to show that $(\bar \xx_\alpha,\bar \yy_\alpha,\bar \zz_\alpha,\bar \rr_\alpha,\bar \sss_\alpha)$ with $(\bar \rr_\alpha,\bar \sss_\alpha)\in\Argmin\pi(\bar \yy_\alpha,\bar \zz_\alpha)$ is also a local optimal solution of  problem \eqref{Penalized_Pb1}. We suppose by contradiction that it is not true. Therefore, there exists a sequence of feasible points $(\xx^n_\alpha,\yy^n_\alpha,\zz^n_\alpha,\rr^n_\alpha,\sss^n_\alpha)_n$ for \eqref{Penalized_Pb1} converging to the point $(\bar \xx_\alpha,\bar \yy_\alpha,\bar \zz_\alpha,\bar \rr_\alpha,\bar \sss_\alpha)$ such that
\[F(\xx^n_\alpha,\yy^n_\alpha)+\alpha((\rr^n_\alpha)^{\top}\zz^n_\alpha+(\sss^n_\alpha)^{\top}(\bb-\A\yy^n_\alpha))<F(\bar \xx_\alpha,\bar \yy_\alpha)+\alpha (\bar \rr_\alpha^{\top}\bar \zz_\alpha+\bar \sss_\alpha^{\top}(\bb-\A\bar \yy_\alpha)),\;\forall n\in\N.\]
It is clear that for all $n\in \N$, $(\xx^n_\alpha,\yy^n_\alpha,\zz^n_\alpha)$  is feasible for \eqref{Penalized_Pb} and we get
\begin{align*}
F(\xx^n_\alpha,\yy^n_\alpha)+\alpha\pi(\yy^n_\alpha,\zz^n_\alpha)&\leq F(\xx^n_\alpha,\yy^n_\alpha)+\alpha((\rr^n_\alpha)^{\top}\zz^n_\alpha+(\sss^n_\alpha)^{\top}(\bb-\A\yy^n_\alpha))\\
&< F(\bar \xx_\alpha,\bar \yy_\alpha)+\alpha (\bar \rr_\alpha^{\top}\bar \zz_\alpha+\bar \sss_\alpha^{\top}(\bb-\A\bar \yy_\alpha))\\
&=F(\bar \xx_\alpha,\bar \yy_\alpha)+\alpha \pi(\bar \yy_\alpha,\bar \zz_\alpha).
\end{align*}
Then the point  $(\bar \xx_\alpha,\bar \yy_\alpha,\bar \zz_\alpha)$ does not  solves \eqref{Penalized_Pb} locally and this is absurd.

As for (ii), let $(\bar \xx_\alpha,\bar \yy_\alpha,\bar \zz_\alpha,\bar \rr_\alpha,\bar \sss_\alpha)$ be a local optimal solution of problem \eqref{Penalized_Pb1}. We show that  $(\bar \xx_\alpha,\bar \yy_\alpha,\bar \zz_\alpha)$ is a local optimal solution of problem \eqref{Penalized_Pb} by making use of the assumption and by observing that from any feasible point $(\xx_\alpha,\yy_\alpha,\zz_\alpha)$ of \eqref{Penalized_Pb}, the point $( \xx_\alpha, \yy_\alpha, \zz_\alpha,\bar \rr_\alpha,\bar \sss_\alpha)$ is feasible for problem \eqref{Penalized_Pb1} and belongs to any neighbourhood of  $(\bar \xx_\alpha,\bar \yy_\alpha,\bar \zz_\alpha,\bar \rr_\alpha,\bar \sss_\alpha)$.
\end{prof}

Combining Propositions \ref{Prop-global-solut-equivalence}, \ref{Prop-local-solut-equivalence}  as well as Proposition \ref{Prop-problem-equivalent} it follows clearly that \eqref{KKT_Reformulation} is globally equivalent to \eqref{Penalized_Pb} if any global optimal solution  $(\bar \xx_\alpha,\bar \yy_\alpha,\bar \zz_\alpha)$ of problem \eqref{Penalized_Pb} belongs to $Z_1^*$. With respect to the local solutions, we obtain the following proposition.
\begin{prop}\label{Prop-local-solut-equivalence1}
 Assume that  for all $\alpha>0$, problem \eqref{Penalized_Pb} possesses an optimal solution which belongs to $Z_1^*$. There exists a scalar $\bar\alpha>0$ such that for all $\alpha\geq\bar\alpha$ such that the following two assertions hold:
\begin{enumerate}
\item[\emph{(i)}] From any local optimal solution  $(\bar \xx_\alpha,\bar \yy_\alpha,\bar \zz_\alpha,\bar \rr_\alpha,\bar \sss_\alpha)$ of \eqref{Penalized_Pb1}, $(\bar \xx_\alpha,\bar \yy_\alpha,\bar \zz_\alpha)$ solves \eqref{KKT_Reformulation} locally if $(\bar \rr_\alpha,\bar \sss_\alpha)\in Argmin\pi( \yy, \zz)$ for all $( \yy, \zz) $ such that $( \bar \xx_\alpha,\yy, \zz)\in Z_1$.
\item[\emph{(ii)}] Suppose that  $(\bar \xx,\bar \yy,\bar \zz)$ is local optimal solution  of \eqref{KKT_Reformulation} for all $\bar \zz\in\R^l $ with $\pi(\bar \yy,\bar\zz)=0$, then $(\bar \xx,\bar \yy,\bar \zz,\bar \rr,\bar \sss)$  with $(\bar \rr,\bar \sss)\in \Argmin \pi(\bar \yy,\bar \zz)$ solves  \eqref{Penalized_Pb1} locally if from any feasible point $(\xx,\yy,\zz)$ of \eqref{Penalized_Pb}, we can find a point $\tilde \zz\in \R^{l}$ such that $(\xx,\yy,\tilde \zz)$ is feasible for \eqref{KKT_Reformulation}.
\end{enumerate}
\end{prop}


To close this section, we would like to discuss some important aspects of the penalization for problem \eqref{Penalized_Pb} considered here. At first, \eqref{Complementarity Reformulation} is not the only way to reformulate the complementarity conditions \eqref{Complementarity}. Various other transformations are possible; see \cite{galantai2012properties} for a large number of such functions and related properties. However, unlike most reformulations that can be found in the latter reference, \eqref{Complementarity Reformulation} enables problem \eqref{Penalized_Pb} to be transformed into the smooth optimization problem \eqref{Penalized_Pb1}.

It is possible to use the bilinear function $(\yy, \zz) \longmapsto \zz^\top(\A\yy-\bb)$ as penalization term in \eqref{Penalized_Pb}, see, e.g., \cite{zemkoho2020theoretical, ye1997exact}. However, its utilization requires the fulfilment of the partial calmness condition, introduced in \cite{ye1995optimality}, which unfortunately does not necessarily hold for \eqref{investigated_Problem}, \eqref{lower-level_problem}.

It is also worth mentioning that exact penalization has been widely used in the context of bilevel programs with lower level problem of the form \eqref{lower-level_problem1}. But it is unclear whether the corresponding results are applicable to \eqref{investigated_Problem}, \eqref{lower-level_problem}. The most recent overview of the topic can be found in \cite{dempeZemkohoBilevelBook2020}.

\section{Semismooth Newton method}\label{Semismooth Newton method}
 Based on the relationship between \eqref{KKT_Reformulation} and \eqref{Penalized_Pb1}, we construct a framework in this section to compute the stationary points of the latter problem using the semismooth Newton method. To begin the process, we write the necessary optimality conditions for \eqref{Penalized_Pb1} by means of the standard Lagrange multipliers rule for smooth optimization problems. As the feasible set of this problem is only described by linear constraints, no constraint qualification is needed.
\begin{thm}\label{Theo_Neccessary_conditions}
Let $(\xx,\yy, \zz,\rr, \sss)\in \R^{2n+3l}$ be a local optimal solution of problem \eqref{Penalized_Pb1} for a fixed value of $\alpha>0$.  Then, there exist $\bm{\lambda}_1\in\R^m_+$, $\bm{\lambda}_6\in\R^n_+$, $\bm{\lambda}_i\in\R^l_+$, $i=2,3,4,5,7$ such that 
\begin{subequations}\label{NOC}
	\begin{align}
		\label{NOC_1}
			&0=\nabla_{\xx} F( \xx,\yy)+\D^\top\bm{\lambda}_1  +\bm{\lambda}_6, \\
		\label{ NOC_2}
&0=\nabla_ {\yy} F( \xx,\yy)-\alpha \A^{\top}\sss +\A^{\top}\bm{\lambda}_2,\\
\label{ NOC_2}
&0=\alpha \rr +\A\bm{\lambda}_6 -\bm{\lambda}_3,\\
\label{ NOC_2}
&0= \alpha  \zz  +\bm{\lambda}_7-\bm{\lambda}_4,\\
\label{ NOC_2}
&0= \alpha (\bb-\A\yy)+\bm{\lambda}_7-\bm{\lambda}_5,\\
\label{ NOC_3}
&0=\A^{\top}\zz+\xx,\\
\label{ NOC_4}
& 0=\rr+\sss-\eee,\\
\label{ NOC_5}
&0=\bm{\lambda}_1^\top(\D\xx-\dd),\;\D\xx\leq \dd,\;\bm{\lambda}_1\geq 0,\\
\label{ NOC_6}
&0=\bm{\lambda}_2^\top(\bb-\A\yy), \bb-\A\yy\geq 0,\;\bm{\lambda}_2\geq 0,\\
\label{ NOC_7}
&0=\bm{\lambda}_3^\top \zz,\;\bm{\lambda}_3\geq 0,\;\zz\geq0,\\
\label{ NOC_6}
&0=\bm{\lambda}_4^\top \rr,\;\bm{\lambda}_4\geq 0,\;\rr\geq0,\\
\label{ NOC_6}
&0=\bm{\lambda}_5^\top \sss,\;\bm{\lambda}_5\geq 0,\;\sss\geq0.
	\end{align}
\end{subequations}
\end{thm}

Next, we transform these conditions into a complete system of equations using the following lemma.  
\begin{lem}[see, e.g., \cite{hintermuller2002primal}]\label{lem1}
For the vectors $\yy_1,\yy_2,\bm{\lambda}\in \R^n$, the system of complementarity conditions $\yy_1\leq \yy_2$, $\bm{\lambda}\geq0$, $\bm{\lambda}^\top(\yy_1-\yy_2)=0$ is equivalent to
$\max\{0,\bm{\lambda}+t(\yy_1-\yy_2)\}-\bm{\lambda}=0$ for any $t>0$.
 \end{lem}
In the process of computing points satisfying \eqref{NOC}, the reformulation of the complementary conditions in Lemma \ref{lem1} has a number of advantages. At first, as it will be clear in the sequel (see, e.g.,  Example \ref{Exa2} and  Theorem \ref{Thm_of_Element_Subdiff}), the possibility to choose any value of $t>0$ provides a level of freedom and flexibility, which can be crucial in solving the resulting system.
 Also, according to \cite{fischer2000merit}, reformulating the complementarity conditions with  linear functions  by means of the Lemma \ref{lem1} can enable the Newton-type method to exhibit the finite termination property while equally demanding slightly weaker assumptions for superlinear convergence. Furthermore, Lemma \ref{lem1} can enable the application of  the Newton method,  to be studied here, to problem in infinite dimensions (e.g., bilevel optimal control problems attracting more and more attention \cite{BenitaMehlitz2016}), through slant  differentiability that generalizes the notion of semismoothness \cite{chen2000smoothing, hintermuller2002primal}.
To solve \eqref{NOC} in the case where the upper level objective function  $F$ is affine linear one, that is, if  $F(\xx,\yy) :=\kk_1^{\top}\xx+\kk_2^{\top}\yy+\kk_3$, $\kk_1, \kk_2\in \R^n$ and $\kk_3\in\R$, Lemma \ref{lem1} allows us to rewrite it as
\begin{equation}\label{SystTo}
\begin{cases}
\B^1X+\B^2\Gamma=\vv,\\
\bm{\Gamma} - \max(0,\bm{\Gamma} +t\bm{\Psi})=0,
\end{cases}
\end{equation}
where $t>0$ is a vector and the max operator is understood componentwise. The given matrices and variables involved in the system \eqref{SystTo} are respectively defined by
\[
\begin{array}{c}
 \B^1:=\left(\begin{array}{ccccc}
\kk_1&\0&\0&\0&\0\\
\kk_2&0&\0&\0&-\alpha \A^{\top}\\
\0&\0&\0&\alpha\I&\0\\
 \0&\0&\alpha\I&\0&\0\\
\0&-\alpha\A&\0&\0&\0\\
\I&\0&\A^{\top}&\0&\0\\
\0&\0&\0&\I&\I
\end{array}\right), \quad
\B^2:=\left(\begin{array}{ccccccc}
\D^\top&\0&\0&\0&\0&\I&\0\\
\0&\A^\top&\0&\0&\0&\0&\0\\
\0&\0&-\I&\0&\0&\A&\0\\
\0&\0&\0&-\I&\0&\0&\I\\
\0&\0&\0&\0&-\I&\0&\I\\
\end{array}\right),\\
\vv:=\left(\begin{array}{c}
\0\\
\0\\
\0\\
\0\\
-\alpha \bb\\
\eee
\end{array}\right), \quad
\bm{\Psi} :=\left(\begin{array}{c}
\D\xx-\dd\\
\A\yy-\bb\\
-\zz\\
-\rr\\
-\sss
\end{array}\right), \quad
X:=\left(\begin{array}{c}
\xx\\
\yy\\
\zz\\
\rr\\
\sss
\end{array}\right),\;\;\mbox{ and }\;\;
\bm{\Gamma}:=\left(\begin{array}{c}
\bm{\lambda}_1\\
\bm{\lambda}_3\\
\bm{\lambda}_4\\
\bm{\lambda}_5\\
\bm{\lambda}_6\\
\bm{\lambda}_7
\end{array}\right).
\end{array}
\]
  If the  matrix $\B^2$ has a full column rank and $(\B^2)^{-1}\B^1$ is a P-matrix, then it can be possible to develop a locally and globally convergent  semismooth Newton scheme to solve \eqref{SystTo}; see, e.g.,  \cite{hintermuller2010semismooth}. However, to expand the number of applications of problem \eqref{investigated_Problem} (cf. next section, for example), we would like the upper level objective function to be a more general twice continuously differentiable function. In this case, solving \eqref{NOC} is equivalent to finding the zeros of the equation
  \begin{equation}\label{Equation Main}
\Phi^{\alpha,t}(\uu):=\begin{pmatrix}
\nabla_{\xx} F(\xx,\yy)+\D^\top\bm{\lambda}_1  +\bm{\lambda}_6\\
\nabla_{\yy} F(\xx,\yy)-\alpha \A^\top\sss + \A^\top\bm{\lambda}_2  \\
\alpha \rr +\A\bm{\lambda}_6-\bm{\lambda}_3\\
\alpha  \zz  +\bm{\lambda}_7-\bm{\lambda}_4\\
\alpha (\bb-\A\yy)+\bm{\lambda}_7-\bm{\lambda}_5\\
\A^{\top}  \zz+\xx\\
\rr+\sss-e\\
\bm{\lambda}_1-\max(0,\bm{\lambda}_1+t_1(\D\xx-\dd))\\
\bm{\lambda}_2-\max(0,\bm{\lambda}_2+t_2(\A\yy-\bb))\\
\bm{\lambda}_3-\max(0,\bm{\lambda}_3+t_3(-\zz))\\
\bm{\lambda}_4-\max(0,\bm{\lambda}_4+t_4(-\rr))\\
\bm{\lambda}_5-\max(0,\bm{\lambda}_5+t_5(-\sss))\\
\end{pmatrix}=0
  \end{equation}
with variable $\uu :=(\xx,\yy,\zz,\rr,\sss,\bm{\lambda}_1,\bm{\lambda}_2,\bm{\lambda}_3,\bm{\lambda}_4,\bm{\lambda}_5,\bm{\lambda}_6,\bm{\lambda}_7)\in \R^{3n+8l+m}$.
It can easily be checked that \eqref{Equation Main} is a square $(3n+8l+m)\times (3n+8l+m)$ system of equations. Additionally, it is well-known that the involved max function is semismooth. Hence, we can solve the system by means of the semismooth Newton method; see, e.g.,  \cite{de1996semismooth, facchinei1997new, qi1993nonsmooth}. In the sequel, we consider the following version of the method due to \cite{de1996semismooth}, as it can be shown that it is globally convergent. 
\begin{algorithm}
\caption{\textsc{: Semismooth Newton algorithm}}
\label{alg1}
\begin{algorithmic}
 \STATE \textbf{Step 0:} Choose $k=0,$ $\uu^0=(\xx^0,\yy^0,\zz^0,\bm{\lambda}_1^0,\bm{\lambda}_2^0,\bm{\lambda}_3^0,\bm{\lambda}_4^0,\bm{\lambda}_5^0,\bm{\lambda}_6^0,\bm{\lambda}_7^0)$, $\delta>0$, $\rho>0$, $p>2$, $\beta\in (0,1/2)$, $\tau=1$, $t>0$,  and $\alpha>0$.
\STATE \textbf{Step 1:} Compute $\Phi^{\alpha,t}(\uu^k)$. If $\tnorm{\Phi^{\alpha,t}(\uu^k)}\leq \delta$, stop.
\STATE \textbf{Step 2:} Select an element $\C^k\in\partial \Phi^{\alpha,t}(\uu^k)$ such that
\begin{equation}\label{Eq13}
\C^k\dd^k+\Phi^{\alpha,t}(\uu^k)=0.
\end{equation}
If \eqref{Eq13} is not solvable or if the condition
\[
\nabla\Psi^{\alpha,t}(\uu^k)^\top \dd^k\leq -\rho\tnorm{\dd^k}^p,
\]
with $\Psi^{\alpha,t}(\uu):=\frac{1}{2}\tnorm{\Phi^{\alpha,t}(\uu)}^2$, is not satisfied, set $\dd^k=-\nabla \Psi^{\alpha,t}(\uu^k)$.
\STATE \textbf{Step 3:} While
\[
\Psi^{\alpha,t}(\uu^k-\tau\nabla\Psi(\uu^k))>\Psi^{\alpha,t}(\uu^k)+\alpha\tau_k\nabla\Psi^{\alpha,t}(\uu^k)^\top\dd^k,
\]
set $\tau_{k+1}=\beta\tau_k$.
\STATE \textbf{Step 4:} Update $ \uu^{k}=\uu^k+\tau_k\dd^k$.
\STATE \textbf{Step 5:} Set $k=k+1$ and go to Step 1.
\end{algorithmic}
\end{algorithm}

\begin{thm}\label{Convergence result}
Fix $\alpha >0$ and $t>0$ and consider a point $\bar \uu$ such that $\Phi^{\alpha,t}(\bar\uu)=\0$. Assume that the upper level objective function $F$ is continuously differentiable and $\nabla F$ is semismooth.
Furthermore, if all matrices $\C\in \partial \Phi^{\alpha,t}(\bar \uu)$ are nonsingular, 
then every sequence generated by Algorithm \ref{alg1} is superlinearly convergent to $\bar \uu$.
Moreover, if  $\nabla F$  is strongly semismooth,
the convergence rate is quadratic.
 \end{thm}
 \begin{prof}
 Follows from \cite{de1996semismooth} by noting that the structure of $\Phi^{\alpha,t}$ ensures that it is semismooth if the upper level objective function $F$ is continuously differentiable and $\nabla F$ is semismooth. Furthermore, $\Phi^{\alpha,t}$ is strongly semismooth provided the same holds for $\nabla F$.
 \end{prof}

In the next lines, we focus on the derivation of sufficient conditions enabling  the nonsingularity of all elements of the generalized Jacobian of  $\Phi^{\alpha,t}$.

\begin{thm} \label{Thm subdiff}
Fixing $\alpha >0$ and $t>0$, the function $\Phi^{\alpha,t}$ \eqref{Equation Main} is strongly semismooth at any point $\uu$ and any matrix from the generalized Jacobian $\partial \Phi^{\alpha,t}(\uu)$ can be written as
\[
\begin{pmatrix}\
\nabla^2_{\xx\xx} F(\xx,\yy)&\nabla^2_{\yy\xx} F( \xx,\yy)&\0&\0&\0&\D^{\top}&\0&\0&\0&\0&\I&\0\\
\nabla^2_{\yy\xx} F( \xx,\yy)&\nabla^2_{\yy\yy} F( \xx,\yy) &\0&\0 &-\alpha \A^{\top}&\0&\A^{\top}&\0&\0&\0&\0&\0\\
\0&\0&\0&\alpha \I&\0&\0&\0&-\I&\0&\0&\A&\0\\
\0&\0&\alpha \I&\0&\0&\0&\0&\0&-\I&\0&\0&\I\\
\0&-\alpha \A&\0&\0&\0&\0&\0&\0&\0&-\I&\0&\I\\
\I&\0&\A^{\top}&\0&\0&\0&\0&\0&\0&\0&\0&\0\\
\0&\0&\0&\I&\I&\0&\0&\0&\0&\0&\0&\0\\
-t_1\pp_1\cdot\D&\0&\0&\0&\0&q_1\cdot\I&\0&\0&\0&\0&\0&\0\\
\0&-t_2\pp_2\cdot\A&\0&\0&\0&\0&q_2\cdot\I&\0&\0&\0&\0&\0\\
\0&\0&t_3\pp_3\cdot\I&\0&\0&\0&\0&q_3\cdot\I&\0&\0&\0&\0\\
\0&\0&\0&t_4\pp_4\cdot\I&\0&\0&\0&\0&q_4\cdot\I&\0&\0&\0\\
\0&\0&\0&\0&t_5\pp_5\cdot\I&\0&\0&\0&\0&q_5\cdot\I&\0&\0
\end{pmatrix}
\]
with the operation $\cdot$ understood as componentwise multiplication and
%
%
\[
\begin{array}{lll}
\pp_1\in \partial\max(0,\bm{\lambda}_1+t_1(\D \xx-\dd)), & \pp_2\in \partial\max(0,\bm{\lambda}_2-t_2(\A\yy-\bb),& \pp_3\in \partial\max(0,\bm{\lambda}_3-t_3\zz),\\
\pp_4\in \partial\max(0,\bm{\lambda}_4-t_4\rr), & \pp_5\in \partial\max(0,\bm{\lambda}_5-t_5\sss), & q_i:=\eee-\pp_i, i\in \{1,\ldots,5\}.
\end{array}
\]
\end{thm}
\begin{prof}
It suffices to calculate the generalized derivative of the last five components of $\Phi^{\alpha,t}$ \eqref{Equation Main}. To illustrate how this can be done, denote by $ Q(\uu):=\bm{\lambda}_1-\max(0,\bm{\lambda}_1+t_1(\D\xx-\dd))$ and observe that
    \[Q(\uu)=\begin{cases}\bm{\lambda}_1& \text{if}\quad \bm{\lambda}_1+t_1(\D\xx-\dd)\leq 0,\\
  -t_1(\D\xx-\dd)&\text{if} \quad \bm{\lambda}_1+t_1(\D\xx-\dd)>0.\\
  \end{cases}\]
 We see that the function $Q$ is differentiable when  $\bm{\lambda}_1+t_1(\D\xx-\dd)>0$ or when $\bm{\lambda}_1+t_1(\D\xx-\dd)<0$. However it is not differentiable when $\bm{\lambda}_1+t_1(\D\xx-\dd)=0$. Hence,
\begin{align*}
\partial Q(\uu)&=\begin{cases}\NN& \;\text{if} \quad \bm{\lambda}_1+t_1(\D\xx-\dd)>0,\\
  \M&\;\text{if} \quad \bm{\lambda}_1+t_1(\D\xx-\dd)<0,\\
  \text{conv}(\{\M,\NN\})&\; \text{if} \quad \bm{\lambda}_1+t_1(\D\xx-\dd)=0,
  \end{cases}\\
  &= \left\{ \pp_1\cdot \NN+(\eee-\pp_1)\cdot \M\;\,\middle|\,\begin{aligned}
   \pp_1=\eee &\; \text{if}\;\bm{\lambda}_1+t_1(\D\xx-\dd)>0\\
  \pp_1=\0 & \;\text{if}\;\bm{\lambda}_1+t_1(\D\xx-\dd)<0\\
  \pp_1\in[\0,\eee] &\; \text{if}\;\bm{\lambda}_1+t_1(\D\xx-\dd)=0
   \end{aligned}
   \right\}
\end{align*}
with $\M$ and $\NN$ respectively defined by
\begin{equation}\label{Eq14}
\begin{array}{rl}
 \M&:=\begin{pmatrix}
\0&\0&\0&\0&\0&\I&\0&\0&\0&\0&\0&\0
\end{pmatrix},\\
\NN&:=\begin{pmatrix}
-t_1\D&\0&\0&\0&\0&\0&\0&\0&\0&\0&\0&\0
\end{pmatrix}
\end{array}
\end{equation}
and $\pp_1\cdot \NN$ corresponds to the matrix whose rows are obtained from the multiplication of the i$_{\text{th}}$ row of $\NN$ and the vector $\pp_1$. The same holds similarly for  $(\eee-\pp_1)\cdot \M$. 
\end{prof}

Next, we identify some scenarios, where the nonsingularity assumption on $\partial \Phi^{\alpha,t}$ required in the convergence Theorem \ref{Convergence result} is satisfied. To proceed, we consider the functions
\[
\begin{array}{l}
h_1(\xx,\yy,\zz,\rr,\sss):=\D\xx-\dd, \;\; h_2(\xx,\yy,\zz,\rr,\sss):=\A\yy-\bb, \;\; h_3(\xx,\yy,\zz,\rr,\sss):=-\zz,\\
h_4(\xx,\yy,\zz,\rr,\sss):=-\rr, \;\; h_5(\xx,\yy,\zz,\rr,\sss):=-\sss,
\end{array}
\]
and the following sets, for $i=1,\ldots,5$:
\[
\left.
\begin{array}{lr}
P_i&:=\tmenge{j\in\{1,\ldots,k_i\}}{0\leq( \bm{\lambda}_i+t_ih_i(\xx,\yy,\zz,\rr,\sss))_j}\\
Q_i&:=\tmenge{j\in\{1,\ldots,k_i\}}{0\geq(\bm{\lambda}_i+t_ih_i(\xx,\yy,\zz,\rr,\sss))_j}
\end{array}
\right\} \mbox{ with } \left\{\begin{array}{l}
                                k_1=m,\\
                                k_2=k_3=k_4=k_5=l.
                              \end{array}
\right.
\]
\begin{thm}\label{theo_regularity_0}
Let $\bar \uu$ be a point that satisfies the optimality conditions in Theorem \ref{Theo_Neccessary_conditions}. Suppose that the matrix $\A$ (with $l=n$) is invertible and $\nabla^2_{xx}F(\bar\xx, \bar\yy)$ has full column rank. Furthermore, assume that the sets  $P_1$, $P_3$, $P_5$, $Q_2$, and $Q_4$ are empty.
Then, all the elements  of $\partial \Phi^{\alpha,t}(\bar \uu)$  are nonsingular.
\end{thm}
\begin{prof}
 Choosing an arbitrary matrix $\C\in \partial \Phi^{\alpha,t}(\bar \uu)$  and consider the homogeneous linear system $\C  \dd= 0$, with a suitably partitioned vector $\dd :=\left(\dd_i\right)_{i=1}^{12}$, we have
 \begin{align}
 \nabla^2_{\xx\xx}F(\bar\xx, \bar\yy)\dd_1+\nabla^2_{\yy\xx}F(\bar\xx, \bar\yy)\dd_2+\D^{\top}\dd_6+\dd_{11}&=\0,\label{Eq 1_0}\\
 \nabla^2_{\xx\yy}F(\bar\xx, \bar\yy)\dd_1+\nabla^2_{\yy\yy}F(\bar\xx, \bar\yy)\dd_2-\alpha \A^{\top}\dd_5+\A^{\top}\dd_7&=\0,\label{Eq 2_0}\\
 \alpha \dd_4+\A\dd_{11}-\dd_8&=\0,\label{Eq 3_0}\\
 \alpha \dd_3+\dd_{12}-\dd_9&=\0,\label{Eq 4_0}\\
 -\alpha \A \dd_2+\dd_{12}-\dd_{10}&=\0,\label{Eq 5_0}\\
 \dd_1+\A^{\top}\dd_3&=0,\label{Eq 6_0}\\
 \dd_4+\dd_5&=\0,\label{Eq 7_0}\\
 -t_1\pp_1\cdot \D\dd_1+(\eee-\pp_1)\dd_6&=\0,\label{Eq 8_0}\\
 -t_2\pp_2\cdot  \A\dd_2+(\eee-\pp_2)\dd_7&=\0,\label{Eq 9_0}\\
  t_3\pp_3\dd_3+(\eee-\pp_3)\dd_8&=\0,\label{Eq 10_0}\\
 t_4\pp_4\dd_4+(\eee-\pp_4)\dd_9&=\0,\label{Eq 11_0}\\
  t_5\pp_5\dd_5+(\eee-\pp_5)\dd_{10}&=\0.\label{Eq 12_0}
 \end{align}
 Since the sets  $P_1,\;P_3,\;P_5,\;Q_2,\;Q_4$ are empty, we get from the definitions of $\pp_1$, $\pp_3$, $\pp_5$, $\pp_2$, and $\pp_4$ that $\pp_1=\pp_3=\pp_5=\0$ and $\pp_2=\pp_4=\eee$.
 It follows then  from \eqref{Eq 7_0}--\eqref{Eq 12_0} and  the invertibility of $\A$ that $\dd_6=\0=\dd_2=\dd_8=\dd_4=\dd_5=\dd_{10}$. Inserting these values in \eqref{Eq 3_0} and in \eqref{Eq 1_0} leads to $\dd_{11}=\0=\dd_1$ since the matrix $ \nabla^2_{\xx\xx}F(\bar\xx, \bar\yy)$ has full column rank.  The reuse of the invertibility of $\A$ in \eqref{Eq 2_0} and \eqref{Eq 6_0} gives $\dd_7=\0=\dd_3$. The equations \eqref{Eq 4_0}--\eqref{Eq 5_0} lead to $\dd_{12}=\0=\dd_9$.
\end{prof}

\begin{exa}
Consider the bilevel optimization problem
 \begin{equation} \label{example_Problem_0}
\begin{array}{rl}
   \underset{x,y}\min & 10x-3x^2+10xy-3y^2\\
 \mbox{ s.t. }& 1\leq x\leq \frac{5}{3}, \;\; y \in \Psi(x)= \underset{y}\Argmin\tmenge{xy}{y\geq 0}.
\end{array}
\end{equation}
The optimal solution of problem $\eqref{example_Problem_0}$ is $(\bar x,\bar y)=(\frac{5}{3},0)$. The objective function of \eqref{example_Problem_0} is concave and the other assumptions of Proposition \ref{Prop-global-solut-equivalence} and Remark \ref{Remark} are satisfied. Therefore there exists
 a scalar $\bar\alpha>0$ such that for all $\alpha\geq\bar\alpha$ 
  we can find $(\bar r,\bar s)\in T$  such that the point $(\frac{5}{3},0,\frac{5}{3},\bar r,\bar s)$ solves the problem
 \begin{equation} \label{example_Problem1_0}
\begin{array}{rl}
   \underset{x,y,z,r,s}\min & 10x-3x^2+10xy-3y^2+\alpha (rz+sy)\\
 \mbox{ s.t. }& 1\leq x\leq \frac{5}{3}, \;\; y\geq 0,\;\; x=z\geq0, \;\; r+s=1,\;\; r\geq0,\;\; s\geq 0.
\end{array}
\end{equation}
The optimal solution of \eqref{example_Problem1_0} is $(\frac{5}{3},0,\frac{5}{3},0,1)$. The corresponding calculations show that the point
\[
\bar\uu=(\bar x,\bar y,\bar z ,\bar r,\bar s,\bar {\bm{\lambda}_1},\bar \lambda_2,\bar \lambda_3,\bar\lambda_4,\bar\lambda_5,\bar\lambda_6,\bar\lambda_7)=\left(\frac{5}{3},0,\frac{5}{3},0,1,0,0,\frac{50}{3}+\alpha,0,\frac{5}{3}\alpha,0,0,0\right)
\]
solves the optimality conditions \eqref{NOC}. Furthermore, it is clear that the matrices $\A=-1$ and $\nabla^2_{xx} F(\bar\xx, \bar\yy)$ are invertible, and the following implications hold
\[
\begin{array}{rllrll}
\bm{\lambda}_1+t_1(\D \bar x-\dd)<0 & \Longrightarrow & \pp_1=\0, & \bar\lambda_3-t_3\bar z<0 &\Longrightarrow& \pp_3=0,\\
\bar\lambda_4-t_4\bar r>0 &\Longrightarrow& \pp_4=1, & \bar\lambda_5-t_5\bar s<0 &\Longrightarrow& \pp_5=0.
\end{array}
\]
Therefore, based on Theorem \ref{theo_regularity_0}, all elements of  $\partial\Phi^{\alpha,t}(\bar\uu)$ are nonsingular for any value of $\alpha\geq \bar \alpha$.
\end{exa}

Obviously, requiring the invertibility of $\A$ in Theorem \ref{theo_regularity_0} is a very strong assumption, as it can enforce the lower level to have a unique solution or feasible point.
 In the sequel, we present a scenario where we avoid to impose this assumption.
\begin{thm}\label{theo_regularity}
Let $\bar \uu$ be a point satisfying the optimality conditions given in Theorem \ref{Theo_Neccessary_conditions}. Suppose that the matrix $\nabla^2_{yy}F(\bar\xx, \bar\yy)$ has full column rank. Furthermore, we assume that the followings sets  $P_1$, $P_2$, $P_5$, $Q_3$, and $Q_4$ are empty.
Then, all the elements  of $\partial \Phi^{\alpha,t}(\bar \uu)$  are nonsingular.
\end{thm}
\begin{prof}
 Proceeding as in the proof of Theorem \ref{theo_regularity_0}, it follows from the emptiness of the sets $P_1,\;P_2,\;P_5,\;Q_3,\;Q_4$ and \eqref{Eq 7_0}--\eqref{Eq 12_0} that $\dd_6=\0=\dd_7=\dd_3=\dd_4=\dd_5=\dd_{10}$ because $\pp_1=\pp_2=\pp_5=\0$ and $\pp_3=\pp_4=\eee$. Inserting these values in \eqref{Eq 2_0} and \eqref{Eq 6_0} leads to $\dd_{1}=\0 =\dd_2$ given that the matrix $ \nabla^2_{yy}F(\bar \xx, \bar\yy)$ is supposed to have full column rank. Inserting again these values in \eqref{Eq 1_0}, \eqref{Eq 3_0}--\eqref{Eq 5_0}  imply that $\dd_{11}=\0=\dd_8=\dd_9=\dd_{12}$.
\end{prof}
\begin{exa}\label{Exa2}
Consider the bilevel optimization problem
 \begin{equation} \label{example_Problem}
\begin{array}{rl}
   \underset{x,y}\min & -3x^2+10xy-3y^2\\
 \mbox{ s.t. }& 1\leq x\leq 2, \;\; y \in \Psi(x)= \underset{y}\Argmin\tmenge{xy}{2\geq y\geq 0}.
\end{array}
\end{equation}
The optimal solution of problem $\eqref{example_Problem}$ is $(\bar x,\bar y)=(2,0)$. The upper level objective function of \eqref{example_Problem} is concave and the other assumptions of Proposition \ref{Prop-global-solut-equivalence} and Remark \ref{Remark} are satisfied. Therefore, there exists
 a scalar $\bar\alpha>0$ such that for all $\alpha\geq\bar\alpha$, 
  we can find $(\bar r,\bar s)=(1/2,2/3,1/2,1/3)$, $\bm{\lambda}_1=(0,0)$, $\bm{\lambda}_2=(\frac{\alpha}{2}-20,0)$, $\bm{\lambda}_3=(\frac{\alpha}{2}+12,\frac{2\alpha}{3}-12)$, $\bm{\lambda}_4=(2\alpha,2/3)$, $\bm{\lambda}_5=(2\alpha,2/3)$, $\lambda_6=12$, $\bm{\lambda}_7=(0,2/3)$ such that the point $\bar u=(\bar x,\bar y,\bar z,\bar r,\bar s,\bm{\lambda}_1,\bm{\lambda}_2,\bm{\lambda}_3,\bm{\lambda}_4,\bm{\lambda}_5,\lambda_6,\bm{\lambda}_7)$ with $\bar z=(2,0)$ solves the optimality conditions  \eqref{NOC}.  Furthermore, it is clear that the matrix $\nabla^2_{yy} F(x,y)=-6$ is full column rank, and since we able to choose $t\in \R_+^5$ with $t_i\neq 0,i=1,\ldots,5$ as desired, the following implications hold
\[
\begin{array}{rllrll}
\bm{\lambda}_1+t_1(\D \bar x-\dd)<0 & \Longrightarrow & \pp_1=\0, & \bm{\lambda}_3-t_3\bar z<0 &\Longrightarrow& \pp_3=\eee,\\
\bm{\lambda}_2+t_2(\A \bar y-\bb)>0 &\Longrightarrow& \pp_2=\0, & \bm{\lambda}_4-t_4\bar r>0 &\Longrightarrow& \pp_4=\eee, \;\; \bm{\lambda}_5-t_5\bar s<0 \Longrightarrow \pp_5=\0.
\end{array}
\]
Therefore, all the assumptions of Theorem \ref{theo_regularity} are fulfilled for $\bar \uu$ and any value of $\alpha\geq \bar \alpha$.
\end{exa}

\begin{rem} It is important to note that the assumptions on  the emptiness of $P^{I_i}_1$ for $j\in\{1,\ldots,5\}$ are not too strong. In fact, they reflect some situations, where the solution w.r.t. the corresponding constraints is nondegenerate or not. For example,  $P^{I_1}_1=\emptyset$ means that  $\bm{\lambda_i}+\ttt_i(\D\xx-\dd)<0$. Clearly, if $\D\xx-\dd<0$, we get the emptiness of $P^{I_1}_1$ from the complementarity conditions.
\end{rem}
 \begin{rem}
As mentioned in the introductory part, the method developed in this section to solve \eqref{investigated_Problem}, \eqref{lower-level_problem} can be applied to \eqref{investigated_Problem}, \eqref{lower-level_problem1}. In fact, for the latter problem, the counterpart of the KKT reformulation \eqref{KKT_Reformulation} can be obtained similarly, without any additional assumption. Hence, we can deduce the optimality conditions and derive a similar generalized derivative as in Theorem \ref{Thm subdiff} and infer  a solution algorithm of  \eqref{investigated_Problem}, \eqref{lower-level_problem1}  while following the same steps.
\end{rem}
\section{Application to the toll-setting problem} \label{sect. app}
We consider the bilevel optimization formulation of the toll-setting problem in transportation. Our aim is to show how the theory discussed in the previous section can be used to solve the problem. First, let us provide a brief description of the problem; for more details, interested readers are referred to \cite{brotcorne2000bilevel,dewez2008new,labbe1998bilevel}, for example. In this problem, the upper level decision-maker corresponds to a road authority or the owner of a highway system allowed to set
tolls on a subset of the links of the network. As for the lower level decision-maker, it is represented by the collection of network users minimizing their travel cost. It is assumed that for a toll selection from the road authority, the network users behave selfishly by trying to minimize their own travel cost. Therefore, the road authority desiring to maximize his/her revenues from tolls will solve the  bilevel program
\begin{equation}\label{Appli1}
\begin{array}{rl}
   \underset{T}\max & \sum_{a\in\mathcal{A}_1} T_ay_a\\[0.5cm]
 \mbox{ s.t. }& \forall a\in \mathcal{A}_1: \;\; T_a\geq l_a,\\[0.5cm]
                         & \underset{y}\min\;\; \sum_{a\in\mathcal{A}_1}(c_a+T_a)y_a+\sum_{a\in\mathcal{A}_2}c_a y_a\\
                         & \forall i\in\mathcal{N},\;\; \forall(j,k)\in \mathcal{O}_\mathcal{D}: \;\; \sum_{a\in i^{+}}y_a^{jk}-\sum_{a\in i^{-}}y_a^{jk}=\begin{cases}1&\text{if} \;i=k_1,\\
                         -1&\text{if}\; i=k_2,\\
                         0&\text{otherwise}, \\
                         \end{cases}\\[0.5cm]
                         &\forall a\in \mathcal{A}: \;\; y_a=\sum_{(j,k)\in \mathcal{O}_\mathcal{D}}d^{jk}y^{jk}_a,\\[0.5cm]
                         &\forall a\in \mathcal{A}, \;\; \forall (j,k)\in \mathcal{O}_\mathcal{D}: \quad y_a^{jk}\geq 0,

\end{array}
\end{equation}

where the involved data and variables are defined as follows:
\begin{longtable}{p{.050\textwidth}  p{.80\textwidth} }
 $\mathcal{A}_1$:& subset of the tolled links,\\
  $\mathcal{A}$:& set of the links,\\
  $\mathcal{A}_2$:& $\mathcal{A}\setminus\mathcal{A}_1$,\\
$T_a$:& toll for link $a\in \mathcal{A}_1$,\\
$c_a$: & (fixed) travel cost for link $a\,(a\in \mathcal{A})$, exclusive of toll,\\
$\mathcal{N}$:& set of nodes,\\
$i^+$:& set of links exiting from node $i \in \mathcal{N}$,\\
$i^-$: & set of links ending at node $i \in \mathcal{N}$,\\
$\mathcal{O}_\mathcal{D}$: & set of origin-destination node pairs,\\
$y_a^{jk}$: &traffic flow from origin $j$ to destination $k$ on link $a$, $(j,k)\in \mathcal{O}_\mathcal{D}$, $ a\in \mathcal{A}$,\\
$d^{jk}$: & proportion of traffic flow demand between origin $j$ and destination $k$,\\
$l_a$:& lower bound on toll for link $a \,(a\in A_1)$,\\
$y_a$: &traffic flow for link $a\in \mathcal{A}$.
\end{longtable}
In order to  keep the variable pattern used so far for the bilevel programming problem  \eqref{investigated_Problem}, \eqref{lower-level_problem}, we  make the following change of variables:
 \[x_a:=\begin{cases}c_a+T_a &\text{if}\;a\in \mathcal{A}_1,\\
 c_a&\text{if}\;a\in \mathcal{A}_2.
 \end{cases}\]
 Therefore, problem \eqref{Appli1}  becomes
\begin{equation}\label{Appli2}
\begin{array}{rl}
   \underset{x}\max & \sum_{a\in\mathcal{A}_1} (x_a-c_a)y_a\\[0.5cm]
 \mbox{ s.t. }& \forall a\in \mathcal{A}_1:\;\; x_a\geq l_a+c_a,\;\; \forall a\in \mathcal{A}_2:\;\; x_a=c_a,\\[0.5cm]
                         & \underset{y}\min\;\; \sum_{a\in\mathcal{A}}x_ay_a\\
                         &\forall i\in\mathcal{N},\;\; \forall(j,k)\in \mathcal{O}_\mathcal{D}: \;\; \sum_{a\in i^{+}}y_a^{jk}-\sum_{a\in i^{-}}y_a^{jk}=\begin{cases}1&\text{if} \;i=j,\\
                         -1&\text{if}\; i=k,\\
                         0&\text{otherwise}, \\
                         \end{cases}\\
                         &\forall a\in \mathcal{A}: \;\; y_a=\sum_{(j,k)\in \mathcal{O}_\mathcal{D}}d^{jk}y^{jk}_a,\\[0.5cm]
                         &\forall a\in \mathcal{A}, \;\; \forall (j,k)\in \mathcal{O}_\mathcal{D}: \;\; y_a^{jk}\geq 0.

\end{array}
\end{equation}

Following the discussion in the previous section, the optimality conditions resulting from the application of Theorem \ref{Theo_Neccessary_conditions} to problem \eqref{Appli2} can be written in the form \eqref{Equation Main} with
\[
 \frac{\partial F(\xx,\yy)}{\partial x_a}=\begin{cases} y_a & \text{if} \; a\in\mathcal{A}_1,\\
0&\text{otherwise}
\end{cases}  \qquad \frac{\partial F(\xx,\yy)}{\partial y_a}=\begin{cases} x_a-c_a &\text{if} \; a\in\mathcal{A}_1,\\
0&\text{otherwise}.
\end{cases}
\]

As for the feasible set of the lower level problem in \eqref{Appli2}, it can be rewritten as $\A_1\yy=\bb_1,\;\yy\geq 0$ with
\[
\begin{array}{l}
\A_1:=\begin{pmatrix}
A_{j_1k_1} &\0&\0&\ldots&\0&\0\\
\0&A_{j_2k_2}&\0&\ldots&\0&\0\\
\vdots&\ddots&\ddots&\ddots&\0&\0\\
\0&\ldots&\0&A_{j_{M-1}k_{M-1}}&\0&\0\\
\0&\ldots&\0&\0&A_{j_Mk_M}&\0\\
d_1^{j_1k_1}\I&d_2^{j_2k_2}\I&\ldots&d_{M-1}^{j_{M-1}k_{M-1}}\I&d_M^{j_Mk_M}\I&-\I\\
\end{pmatrix}, \; \yy :=\begin{pmatrix} y_a^{jk}\\
y_a\end{pmatrix}_{a,j,k},\; \bb_1  :=\begin{pmatrix}
b^i_{j_sk_s}\\
\0
\end{pmatrix}_{s}\\
\forall s=1,\ldots,M, \;\;\forall (j_s,k_s)\in \mathcal{O}_\mathcal{D},\;\; \forall i\in\mathcal{N}: \quad b^i_{j_sk_s}=\begin{cases}1&\text{if} \;i=j_s,\\
                         -1&\text{if}\; i=k_s,\\
                         0&\text{otherwise}.\\
                         \end{cases}
\end{array}
\]
The matrix $A_{j_sk_s}$, $s=1\ldots M$ represents the incidence matrix of the graph with origin $j_i$ and  destination $k_i$. Here we have set $M:=|\mathcal{O}_\mathcal{D}|$ to be the number of start-destination node pairs. We want to draw the attention of the reader on the fact that, for each individual graph an equation of the traffic flow condition can be omitted, i.e. a row in $\A_1$ belonging to the matrix $A_{jk} $ and the corresponding component in the vector $\bb_1$ can be deleted, since the systems of equations are super-determined for a single flow. Furthermore the number of variables can reduce if the vehicles  from node $j$ to node $ k$ are allowed to use during their journey  only the road which are useful.

To get the lower level feasible set in the format in \eqref{investigated_Problem}, \eqref{lower-level_problem}, the above described system $\A_1\yy=\bb_1,\;\yy\geq 0$ is completely transformed into inequality conditions.
In the same vein, the matrix $\D$ from the upper level constraint is constructed by replacing the equality constraints  by inequalities. It is obtained from the constraints $-x_a\leq -l_a-c_a$ for $a\in \mathcal{A}_1$ and  $x_a\leq c_a$,  $-x_a\leq -c_a$ for $a\in \mathcal{A}_2$. This means that for all $j\in\mathcal{A}$,
\[
\D_{ij} := \begin{cases}-\delta_{ij}&\text{if}\; i\in\mathcal{A}_1,\\
\begin{pmatrix}
\delta_{ij}\\
-\delta_{ij}
\end{pmatrix}&\text{if} \;i\in\mathcal{A}_2, \end{cases}\qquad \dd_{i} := \begin{cases}-l_a-c_a&\text{if}\; i\in\mathcal{A}_1,\\
\begin{pmatrix}
c_a\\
-c_a
\end{pmatrix}&\text{if} \;i\in\mathcal{A}_2. \end{cases}
\]

 It is worth recognizing that transforming an equality into two inequalities is not computationally ideal. Nevertheless, we do so here just to be faithful to the model \eqref{investigated_Problem}, \eqref{lower-level_problem} and the corresponding analysis conducted in the previous sections.


To implement Step 2 of Algorithm \ref{alg1} on the toll-setting problem, we should be able to calculate an element of $\partial \Phi^{\alpha,t}(\uu^k)$. Our first task is therefore to show how this
can be accomplished. To proceed, we define the vectors
\begin{align*}
X^1(\uu)&:= \bm{\lambda}_1+t_1(\D\xx-\dd),& X^3(\uu)&:=\bm{\lambda}_3-t_3\zz, & X^5(\uu)&:=\bm{\lambda}_5-t_5\sss.\\
X^2(\uu)&:=\bm{\lambda}_2+t_2(\A\yy-\bb),&  X^4(\uu)&:=\bm{\lambda}_4-t_4\rr,
\end{align*}
\begin{thm}\label{Thm_of_Element_Subdiff}
 The following matrix  is an element of $\partial \Phi^{\alpha,t}(\uu)$ for $\uu\in \R^{3n+8l+m}$:
\[
H:=\begin{pmatrix}
\0&\nabla^2_{\yy\xx} F( \xx,\yy)&\0&\0&\0&\D^{\top}&\0&\0&\0&\0&\I&\0\\
\nabla^2_{\yy\xx} F( \xx,\yy)&\0&\0&\0 &-\alpha \A^{\top}&\0&\A^{\top}&\0&\0&\0&\0&\0\\
\0&\0&\0&\alpha \I&\0&\0&\0&-\I&\0&\0&\A&\0\\
\0&\0&\alpha \I&\0&\0&\0&\0&\0&-\I&\0&\0&\I\\
\0&-\alpha \A&\0&\0&\0&\0&\0&\0&\0&-\I&\0&\I\\
\I&\0&\A^{\top}&\0&\0&\0&\0&\0&\0&\0&\0&\0\\
\0&\0&\0&\I&\I&\0&\0&\0&\0&\0&\0&\0\\
K^1&\0&\0&\0&\0&\tilde K^1 &\0&\0&\0&\0&\0&\0\\
\0&K^2&\0&\0&\0&\0&\tilde K^2&\0&\0&\0&\0&\0\\
\0&\0&K^3&\0&\0&\0&\0&\tilde K^3&\0&\0&\0&\0\\
\0&\0&\0&K^4&\0&\0&\0&\0&\tilde K^4&\0&\0&\0\\
\0&\0&\0&\0&K^5&\0&\0&\0&\0&\tilde K^5&\0&\0
\end{pmatrix}.
\]
Here, we have
$
\frac{\partial^2 F}{\partial x_a\partial y_a} = \frac{\partial^2 F}{\partial y_a\partial x_a} = \left\{\begin{array}{ll}
                                                                                                         1 & \mbox{if }\; a\in \mathcal{A}_1,\\
                                                                                                         0 & \mbox{otherwise},
                                                                                                       \end{array}\right.
$
\[
\begin{array}{lll}
K^1 &:=&  -\frac{t_1}{2}\left(\D+\sgn(X^1)\cdot \D\right),\\[1ex]
K^2 &:=&  -\frac{t_2}{2}\left(\A+\sgn(X^2)\cdot \A\right),
\end{array}\quad
\begin{array}{lll}
K^i &:=&  -\frac{t_i}{2}\left(\I+\sgn(X^i)\cdot \I\right), \;i=3,4,5,\\[1ex]
\tilde K^i &:=&  -\frac{1}{2}\left(\I-\sgn(X^i)\cdot \I\right),\;i=1,\ldots,5.
\end{array}
\]
\end{thm}
\begin{prof}
Observe that the first $3n+4l$ components of the function $\Phi^{\alpha,t}$ are differentiable. Hence, we focus our attention on the components which contain the $\max$ operator. In order to prove that $H\in\partial \Phi^{\alpha,t}(\uu)$, we will  build a sequence of points $(\tilde\uu^k)_k$ and $(\tilde{\tilde\uu}^k)_k$,  where $ \Phi^{\alpha,t}$ is differentiable  at the points $\uu+\varepsilon^k\tilde \uu$ and $\uu+\varepsilon^k\tilde{\tilde\uu}^k$ and it holds $\nabla \Phi^{\alpha,t} (\uu+\varepsilon^k\tilde \uu)\rightarrow H^1$ and $\nabla \Phi^{\alpha,t} (\uu+\varepsilon^k\tilde{\tilde\uu}^k)\rightarrow H^2$ for $k\rightarrow\infty$. After that  we will show that $H$ is an element of the convex hull of $H^1$ and $H^2$. Let
\[
\Delta :=\left\{i\in\{1,\ldots,3n+8l+m\}\;|\;X_i^j(\uu)=0,\;j\in\{1,\ldots,5\}\right\}.
\]
 We consider some vectors
\[
\tilde\uu:=(\tilde\xx,\tilde\yy,\tilde\zz,\tilde\rr,\tilde\sss,\tilde{\bm{\lambda}}_1,\tilde{\bm{\lambda}}_2,\tilde{\bm{\lambda}}_3,\tilde{\bm{\lambda}}_4,\tilde{\bm{\lambda}}_5,\tilde{\bm{\lambda}}_6,\tilde{\bm{\lambda}}_7)\;
\mbox{ and } \;
\tilde{\tilde\uu}:=(\tilde{\tilde\xx},\tilde{\tilde\yy},\tilde{\tilde\zz},\tilde{\tilde\rr},\tilde{\tilde\sss},\tilde{\tilde{\bm{\lambda}}}_1,\tilde{\tilde{\bm{\lambda}}}_2,\tilde{\tilde{\bm{\lambda}}}_3,\tilde{\tilde{\bm{\lambda}}}_4,\tilde{\tilde{\bm{\lambda}}}_5,\tilde{\tilde{\bm{\lambda}}}_6,\tilde{\tilde{\bm{\lambda}}}_7)
\]
from $\mathbb{R}^{3n+8l+m}$, such that $\tilde \uu_j=0=\tilde{\tilde \uu}_j$ for $j\notin\Delta$ and $\tilde \uu_j=1$, $\tilde{\tilde \uu}_j=-1$ for $j\in \Delta $.

If there exists $i_0\in \{3n+4l+1,\ldots,3n+8l+m\}$ such that $X_{i_0}^j(\uu)\neq 0$, $j=1,\ldots,5$, suppose, without loss of generality that $X_{i_0}^1(\uu)\neq 0$ (because we can apply the same reasoning on all  $j=1,\ldots,5$). Then $\Phi^{\alpha,t}_{i_0}$ is differentiable at $\uu$ and we get from the proof of Theorem \ref{Thm subdiff}, cf. \eqref{Eq14}, that
\[\partial \Phi^{\alpha,t}_{i_0}(\uu)=
\begin{cases}
\NN_{i_0}&\;\text{if}\quad X^1_{i_0}>0,\\
\M_{i_0}&\;\text{if}\quad X^1_{i_0}<0.\\
\end{cases}
\]
This gives the row $H_{i_0}$ since
 $\sgn(X_{i_0}^1)=\pm 1$ and $K^1_{i_0}=-t_1 \D_{i_0}$ if $X_{i_0}^1>0$ and $K^1_{i_0}=0$ if $X_{i_0}^1<0$.

Otherwise, suppose that for all $i\in \{3n+4l+1,\ldots,3n+8l+m\}$,  $X_{i}^j= 0$, $j=1,\ldots,5$. We have for example $X^1=\0$, then $\sgn(X^1)=\0$. We consider the sequences $\tilde\uu_i^k=\tilde\uu_i+\varepsilon^k$ and  $\tilde{\tilde \uu}^k_i=\uu_i+\varepsilon^k$, $k\in\N$, and $\varepsilon^k\downarrow0$.  It holds
\begin{align*}
X^1(\uu+\varepsilon^k\tilde \uu)&=(\bm{\lambda}_1+\varepsilon^k\tilde{\bm{\lambda}_1})+t_1(\D(\xx+\varepsilon^k\tilde{\xx})-\dd)\\
&=\underbrace{\bm{\lambda}_1+t_1(\D\xx-\dd)}_{=0}+\varepsilon^k\tilde{\bm{\lambda}_1}+t_1\varepsilon^k\D\tilde{\xx}\\
&=\varepsilon^k(\tilde{\bm{\lambda}_1}+t_1\D\tilde{\xx})>0\quad \text{for $t_1$ small enough}.
\end{align*}

In the same vein, we have $X^1(\uu+\epsilon^k\tilde{\tilde \uu})=\epsilon^k(\tilde{\tilde{\bm{\lambda}_1}}+t_1\D\tilde{\tilde{\xx}})<0\quad \text{for $t_2$ small enough}$. Consequently, the mapping  $\Phi_1^{\alpha,t}:\bm{\uu}\mapsto\bm{\lambda}_1-\max(0,\bm{\lambda}_1+t_1(\D\xx-\dd))$ is at the points $\uu+\epsilon^k\tilde \uu$ and $\uu+\varepsilon^k\tilde{\tilde \uu}$ differentiable and we get
\begin{align*}
\nabla \Phi_1^{\alpha,t}(\uu+\varepsilon^k\tilde \uu)&\underset{k\to \infty}\to \begin{pmatrix}
\0&\0&\0&\0&\0&\I&\0&\0&\0&\0&\0&\0
\end{pmatrix}:=H^1\\
\nabla \Phi_1^{\alpha,t}(\uu+\varepsilon^k\tilde{\tilde \uu})&\underset{k\to \infty}\to\begin{pmatrix}
-t_1\D&\0&\0&\0&\0&\0&\0&\0&\0&\0&\0&\0
\end{pmatrix}:=H^2
\end{align*}

and clearly the matrix
\[\begin{pmatrix}
\frac{-t_1}{2}\D&\0&\0&\0&\0&\frac{1}{2}\I&\0&\0&\0&\0&\0&\0
\end{pmatrix}\]
belongs to the convex hull of $H^1$ and $H^2$. You see that if you apply the same reasoning on all the others rows you will recover the matrix $H$ given in the theorem.
\end{prof}

The particularity of the matrix $H$ constructed in this theorem is the fact that it is the limit of the gradient of the following differentiable approximation
\begin{equation}\label{Approx of Phi}
\Phi^{\alpha,t}_\varepsilon(\uu):=\begin{pmatrix}
\nabla_{\xx} F(\xx,\yy)+\D^\top\bm{\lambda}_1  +\bm{\lambda}_6\\
\nabla_{\yy} F(\xx,\yy)-\alpha \A^\top\sss + \A^\top\bm{\lambda}_2  \\
\alpha \rr +\A\bm{\lambda}_6-\bm{\lambda}_3\\
\alpha  \zz  +\bm{\lambda}_7-\bm{\lambda}_4\\
\alpha (\bb-\A\yy)+\bm{\lambda}_7-\bm{\lambda}_5\\
\A^{\top}  \zz+\xx\\
\rr+\sss-e\\
\bm{\lambda}_1-\frac{1}{2}(X^1+\sqrt{(X^1)^2+\varepsilon})\\
\bm{\lambda}_2-\frac{1}{2}(X^2+\sqrt{(X^2)^2+\varepsilon})\\
\bm{\lambda}_3-\frac{1}{2}(X^3+\sqrt{(X^3)^2+\varepsilon})\\
\bm{\lambda}_4-\frac{1}{2}(X^4+\sqrt{(X^4)^2+\varepsilon})\\
\bm{\lambda}_5-\frac{1}{2}(X^5+\sqrt{(X^5)^2+\varepsilon})\\
\end{pmatrix}
\end{equation}
of the function $\Phi^{\alpha, t}$. Here, $\varepsilon > 0$ and for $x\in \R$,
\[
x+\sqrt{x^2+\varepsilon} \;\, \longrightarrow \;\, \max (0, x) \;\; \mbox{ as } \;\, \varepsilon \;\, \longrightarrow \;\, 0.
\]
We recall that the square and the square-root applied on vectors in \eqref{Approx of Phi} are understood componentwise. It can easily be checked that the gradient of $\Phi^{\alpha,t}_\varepsilon$ can be written as
\[\nabla\Phi^{\alpha,t}_\varepsilon(\uu)=\begin{pmatrix}
\nabla^2_{\xx\xx} F(\xx,\yy)&\nabla^2_{\yy\xx} F( \xx,\yy)&\0&\0&\0&\D^{\top}&\0&\0&\0&\0&\I&\0\\
\nabla^2_{\yy\xx} F( \xx,\yy)&\nabla^2_{\yy\yy} F( \xx,\yy) &\0&\0 &-\alpha \A^{\top}&\0&\A^{\top}&\0&\0&\0&\0&\0\\
\0&\0&\0&\alpha \I&\0&\0&\0&-\I&\0&\0&\A&\0\\
\0&\0&\alpha \I&\0&\0&\0&\0&\0&-\I&\0&\0&\I\\
\0&-\alpha \A&\0&\0&\0&\0&\0&\0&\0&-\I&\0&\I\\
\I&\0&\A^{\top}&\0&\0&\0&\0&\0&\0&\0&\0&\0\\
\0&\0&\0&\I&\I&\0&\0&\0&\0&\0&\0&\0\\
G^1&\0&\0&\0&\0&\tilde G^1 &\0&\0&\0&\0&\0&\0\\
\0&G^2&\0&\0&\0&\0&\tilde G^2&\0&\0&\0&\0&\0\\
\0&\0&G^3&\0&\0&\0&\0&\tilde G^3&\0&\0&\0&\0\\
\0&\0&\0&G^4&\0&\0&\0&\0&\tilde G^4&\0&\0&\0\\
\0&\0&\0&\0&G^5&\0&\0&\0&\0&\tilde G^5&\0&\0
\end{pmatrix},\]
where
\begin{align*}
G^1&:= \frac{-t_1}{2}(\D+(X^1 \div L^1)\cdot \D),& G^i&:=\frac{-t_i}{2}(\I+(X^i \div L^i)\cdot \I), \;i=3,4,5,\\
G^2&:=\frac{-t_2}{2}(\A+(X^2 \div L^2)\cdot \A),&  \tilde G^i&:=\frac{-1}{2}(\I-(X^i \div L^i)\cdot \I),\;i=1,\ldots,5.
\end{align*}
Here, $L^i_j:=\sqrt{X_j^i+\varepsilon}\quad j=1.\ldots.n,\; i=1.\ldots.5$. The signs $\div$ and $\cdot$ are used to denote the componentwise division and multiplication, respectively. This means we have  two options to evaluate the gradient of the merit function in Step 2 of Algorithm \ref{alg1}: either compute the gradient of  $\nabla \Psi^{\alpha, t}$  directly if possible or the derivative of the approximation function $\Phi^{\alpha,t}_{\varepsilon}$ for small values of $\varepsilon$.

\begin{figure}[H]
\centering
\begin{subfigure}{.5\textwidth}
\centering
\begin{tikzpicture}[scale=1.8]
	\node (1) at (1,1) [circle,draw] {1};
	\node (2) at (3,1) [circle,draw] {2};
	\node (3) at (3,-1) [circle,draw] {3};
	\node (4) at  (1,-1) [circle,draw] {4};
	\node (5) at (-1,0) [circle,draw] {5};
	
	\draw[->, thick] (1) to node[above] {$c_1=2$} (2);
	 \draw[->, thick] (2) to node[above] {$c_4=0$} (3);
	 \draw[->, thick] (3) to node[above] {$c_6=2$}(4);
	 \draw[->, thick](4) to node[above] {$c_8=0$} (5);
	 \draw[->, thick] (1) to node[above] {$c_3=5$} (5);
	 \draw[-, thick] (1) to  (2,0);
	  \draw[->, thick] (2,0) to  node[above] {$c_2=6$} (3);
	 \draw[-, thick] (2) to   (2,0);	
	  \draw[->, thick]  (2,0) to  node[above] {$c_5=4$}(4);
	  \draw[<-, thick] (5) to (-1,-2);
	 \draw[-, thick] (-1,-2) to node[above] {$c_7=6$} (3,-2);
	  \draw[-, thick] (3,-2) to (3);
	
	  	
\end{tikzpicture}
  \caption{Network 1} \label{fig:test1}
\end{subfigure}%
\begin{subfigure}{.5\textwidth}
\centering
\begin{tikzpicture}[scale=1.3]
	\node (1) at (-1,2) [circle,draw] {1};
	\node (2) at (4,2) [circle,draw] {2};
	\node (3) at (0,0) [circle,draw] {3};
	\node (4) at  (3,0) [circle,draw] {4};
	\node (5) at (-1,-2) [circle,draw] {5};
	\node (6) at (4,-2) [circle,draw] {6};
	
	\draw[->, thick] (1) to node[above] {$c_1=8$} (2);
	 \draw[->, thick] (1) to node[above] {$c_2=2$} (3);
	 \draw[->, thick] (5) to node[above] {$c_5=3$}(3);
	 \draw[->, thick](5) to node[above] {$c_7=6$} (6);
	 \draw[->, thick] (4) to node[above] {$c_6=1$} (6);	
	 \draw[->, thick] (3) to node[above] {$c_4=\frac{x_4}{2}$} (4);	
	   \draw[<-, thick] (4) to  node[above] {$c_3=1$} (2); 	 	
\end{tikzpicture}
 \caption{Network 2}\label{fig:test2}
\end{subfigure}
\caption{Networks with 5 and 6 nodes, respectively}\label{fig:test}
\end{figure}
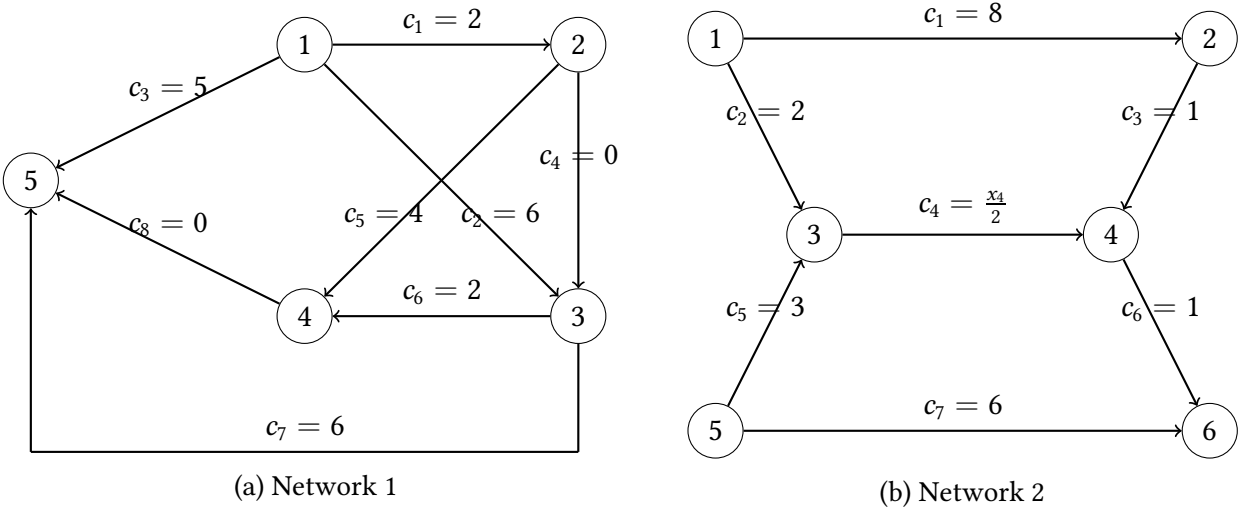

To illustrate the application of Algorithm \ref{alg1} to problem \eqref{Appli2}, we consider the examples of networks in Figure \ref{fig:test1} and \ref{fig:test2} , which are taken from \cite{colson2002bilevel}.
The origin-destination nodes for Network 1 (resp.  Network 2) is $(1,5)$ (resp. $\{(1,2),(5,6)\}$).
The formulation of the toll problem related to these networks respectively leads to bilevel programming problems explicitly written as
\[
\begin{array}{c}
 \begin{array}{rl}
   \underset{\xx,\yy}\max & (x_3-5)y_3+x_4y_4+x_8y_8\\
 \mbox{ s.t. }& \forall a\in \mathcal{A}_2, x_a=c_a,\\
 &\forall a\in \mathcal{A}_1, x_3\geq 5, x_4\geq 0,x_8\geq 0,\\
                         & \underset{y}\min \,x^{\top}y\\
                         &y_1+y_2+y_3=1,\\
                         &y_8+y_7+y_3=1,\\
                         &y_4+y_5-y_1=0,\\
                         &y_6+y_7-y_2-y_4=0,\\
                         &y_8-y_5-y_6=0,\\
                         &y_i\geq0,\;\forall i=1,\ldots,8,

\end{array} \qquad\;\; \mbox{ and } \qquad\;\;
\begin{array}{rl}
   \underset{\xx,\yy}\max &\frac{1}{2} x_4y_4\\
 \mbox{ s.t. }& \forall a\in \mathcal{A}_2, x_a=c_a,\\
                   &x_4=x_8,\\
                          & \underset{y}\min \,x^{\top}y\\
                         &y_1+y_2=1,\\
                         &y_5+y_7=1,\\
                         &y_2=y_3=y_4,\\
                         &y_5=y_6=y_8,\\
                         &y_i\geq0,\;\forall i=1,\ldots,8.

\end{array}
\end{array}
\]
We have $\mathcal{A}_1=\{3,4,8\}$ and  $\mathcal{A}_2=\{1,2,5,6,7\}$ for Network 1 and $\mathcal{A}_1=\{4\}$ and  $\mathcal{A}_2=\{1,2,3,5,6,7\}$ for Network 2. The variables $y_8$ and $y_4$ in the mathematical formulation of the toll problem related to the Network 2 denote $y^{12}_4$ and $y^{56}_4$, respectively and  the travel cost $c_a $ for each link $(a\in \mathcal{A})$
 can be obtained on each corresponding network. Next, we  apply Algorithm \ref{alg1} to solve these two problems. To proceed, we consider the parameters:
\[
t(1)=0.045,\; t(2)=0.049,\; t(3)=0.025,\; t(4)=0.005, \;t(5) =0.0025,\; \varepsilon =0.01, \; \delta = 10^{-6}.
\]
For these examples, we set the algorithm to stop if $\tnorm{\Phi^{\alpha,t}(\uu^k)}\leq \delta$ or the iteration index $k$ reaches 50. As starting point  for Network 1, we choose
\[
\begin{array}{l}
 \xx_0=(2, 6, 5, 5, 4, 2, 6, 5), \;\; \yy_0=( 0, 1, 0, 0, 0, 1, 0, 0), \;\; \zz_0 = \A y_0-b=\bm{\lambda}_2=\bm{\lambda}_3,\\
 \rr_0 =\bm{\lambda}_4=\bm{\lambda}_7=\0_{l\times 1}, \;\; \sss_0 = \eee=\bm{\lambda}_5, \;\; \bm{\lambda}_1 = |\D \xx_0-\dd|, \;\; \bm{\lambda}_6=\0_{n\times 1},
\end{array}
\]
and for $\alpha=0.45$, we obtain $F(\xx^*,\yy^*)=  -7.0346$ and $f(\xx^*,\yy^*)= 12.9927$,
while the known values in the literature are $F(\xx^*,\yy^*)= -7$ and $f(\xx^*,\yy^*)= 13$, see \cite{colson2005trust}.
For Network 2, we choose
\[
\begin{array}{l}
 \xx_0=(8, 2, 1, 0, 3, 1, 6, 0), \;\; \yy_0=\0_{n\times 1}, \;\; \zz_0 = \A y_0-b=\bm{\lambda}_2=\bm{\lambda}_3,\\
 \rr_0 =\bm{\lambda}_4=\bm{\lambda}_7=\0_{l\times 1}, \;\; \sss_0 = \eee=\bm{\lambda}_5, \;\; \bm{\lambda}_6=\0_{n\times 1},
\end{array}
\]
as a starting point and for $\alpha=4.791$, we get $F(\xx^*,\yy^*)=  -11.7003$ and $f(\xx^*,\yy^*)= 34.0099$.
The known values in the literature are $F(\xx^*,\yy^*) =  -11$ and $f(\xx^*,\yy^*)= 34$; cf. latter reference.

{
\begingroup
\setlength\bibitemsep{-1.5pt}
\setlength{\parskip}{-1.5pt}
\printbibliography
\endgroup
}
\end{document}